\documentclass[11pt]{article}

\usepackage{amssymb}

\usepackage{amsmath}

\usepackage{amsfonts}

\usepackage{epsfig}

\usepackage{graphics}

\newtheorem{lemmamr}{Lemma}

\newtheorem{conjecturemr}[lemmamr]{Conjecture}

\begin{document}

\title{Secondary terms in the number of vanishings of quadratic
twists of elliptic curve $L$-functions}

\author{J.\ Brian\ Conrey, Atul Pokharel, Michael\ O.\ Rubinstein \\and Mark Watkins}

\maketitle

\begin{abstract}

We examine the number of vanishings of
quadratic twists of the $L$-function associated to an
elliptic curve. Applying a conjecture for the full asymptotics
of the moments of critical $L$-values we obtain a conjecture for
the first two terms in the ratio of the number of vanishings
of twists sorted according to arithmetic progressions.
\end{abstract}


\section{Introduction}
\label{section:intro}

Let $E$ be an elliptic curve over $\mathbb{Q}$ with associated $L$-function
given by
\begin{align}
\label{eq:L_E}
    L_{E}(s) & =\sum_{n=1}^\infty\frac{a_n}{n^s}
    =\prod_{p\mid \Delta}
    \left(1-a_p p^{-s}\right)^{-1}
    \prod_{p\nmid \Delta}
    \left(1-a_p p^{-s}+p^{1-2s}\right)^{-1}\\
    & =
    \prod_{p} {\mathcal{L}}_p(1/p^s)
    , \quad \quad \Re(s) > 3/2.
\end{align}
Here, $\Delta$ is the discriminant of $E$, and $a_p=p+1-\#E(\mathbb{F}_p)$,
with $\#E(\mathbb{F}_p)$ the number of points, including the point at infinity,
of $E$ over $\mathbb{F}_p$.
$L_E(s)$ has analytic continuation to $\mathbb{C}$ and satisfies a
functional equation~\cite{kn:wiles95}~\cite{kn:taywil95}~\cite{kn:BCDT01} of the form
\begin{equation}
\label{eq:10a}
    \left(\frac{2\pi}{\sqrt{Q}}\right)^{-s}{\rm
    \Gamma}(s)L_E(s)=w_E\left(\frac{2\pi}{\sqrt{Q}}\right)^{s-2}{\rm\Gamma}(2-s)L_E(2-s),
\end{equation}
where $Q$ is the conductor of the elliptic curve $E$ and $w_E=\pm 1$.

Let
\begin{equation}
    \label{eq:functional eqn}
    L_{E}(s,\chi_d)=\sum_{n=1}^{\infty}\frac{a_n\chi_d(n)}{n^s}
\end{equation}
be the $L$-function of the elliptic curve $E_d$, the quadratic
twist of $E$ by the fundamental discriminant $d$.
If ($d,Q)=1$, then $L_{E}(s,\chi_d)$ satisfies the functional equation
\begin{equation}
    \label{eq:twisted functional eqn}
    \left(\frac{2\pi}{\sqrt{Q}|d|}\right)^{-s}
    {\rm \Gamma}(s)L_{E}(s,\chi_d)=\chi_d(-Q)w_E\left(\frac{2\pi}{\sqrt{Q}|d|}\right)^{s-2}
    {\rm \Gamma}(2-s)L_{E}(2-s,\chi_d).
\end{equation}

In \cite{kn:ckrs00} and \cite{kn:ckrs04} conjectures, modeled after corresponding
theorems in random matrix theory, are stated concerning the distribution of values of
$L_{E}(1,\chi_d)$ with an application made to counting the number of
vanishings of $L_{E}(1,\chi_d)$.
We focus on the case $w_E \chi_d(-Q)=1$, since
otherwise $L_{E}(1,\chi_d)$ is trivially equal to zero.
One quantity studied concerns the ratio of the number of
vanishings sorted according to residue classes mod $q$ for a fixed prime $q \nmid Q$.
Let
\begin{equation}
    \label{eq:Rq}
    R_q(X)=
    \frac{
        \sum_{{|d|<X, w_E \chi_d(-Q)=1  \atop L_E(1,\chi_d)=0} \atop \chi_d(q)=1} 1
    }
    {
        \sum_{{|d|<X, w_E \chi_d(-Q)=1  \atop L_E(1,\chi_d)=0} \atop \chi_d(q)=-1} 1
    }
\end{equation}
be the ratio of the number of vanishings of $L_E(1,\chi_d)$ sorted according
to whether $\chi_d(q)=1$ or $-1$.

By looking at this ratio, certain elusive and mysterious
quantities that appear in the asymptotics for both the numerator
and denominator cancel each other out and one is left with a
precise prediction for its limit. Let
\begin{equation}
    R_q = \left( \frac{q+1-a_q}{q+1+a_q} \right)^{1/2}.
\end{equation}
A conjecture from~\cite{kn:ckrs00} asserts that, for $q \nmid Q$,
\begin{equation}
    \lim_{X\rightarrow \infty}
    R_q(X)=R_q.
\end{equation}
It is believed that this continues to hold if the set
of quadratic twists is restricted to subsets such as
$d<0$ or $d>0$, or to $|d|$ prime,
though in the latter case we must be sure to rule out there
being no vanishings at all due to arithmetic reasons~\cite{kn:delaunay04}.

Numerical evidence for three elliptic curves is presented in \cite{kn:ckrs00}
and confirms this prediction. However, even taking $X$ of size roughly $10^9$
(and, in that paper, $d<0$ and $|d|$ prime), the numeric value of the
ratio was found in that paper to agree with the predicted value to about two
decimal places. In other cases, when $a_q$ of $L_E(S)$ in~(\ref{eq:L_E}) equals 0,
the numeric value of $R_q(X)$ compared to the predicted limit $R_q$
to three or more decimal places.

In this paper we examine secondary terms in the above conjecture
applying new conjectures~\cite{kn:cfkrs} for the full asymptotics of the
moments of $L_E(1,\chi_d)$. We obtain a conjectural formula for the next to leading
term in the asymptotics for $R_q(X)$. It is of size $O(1/\log(X))$ and explains the
slow convergence to the limit $R_q$. We also explain in Section 3 the tighter fit
when $a_q=0$.

While the main term, $R_q$, in the above conjecture is robust and does not
depend heavily on the set of $d$'s considered, the secondary terms are
more sensitive, for example, to the residue classes of $d$ modulo the primes
that divide $Q$.
Therefore, for simplicity we focus on the following
dense collection of fundamental discriminants $d$. Assume that $Q$ is squarefree and let
\begin{equation}
\label{eq:S-}
    S^-(X) = S^-_E(X) = \{ -X \leq d < 0; \chi_d(p) =  -a_p \ \ \text{for all $p \mid Q$}\}
\end{equation}
For curves of prime conductor $Q$ we also consider the set of fundamental discriminants
\begin{equation}
\label{eq:S+}
    S^+(X) = S^+_E(X) = \{ 0 < d \leq X; \chi_d(Q) =  a_Q \}.
\end{equation}
These sets of discriminants are also chosen because they allow us to efficiently compute
$L_E(1,\chi_d)$ using a relationship to the coefficients of certain modular forms of weight
$3/2$ that has been worked out explicitly for many examples by Tornaria and
Rodiguez-Villegas~\cite{kn:rt04} (see~\cite{kn:ckrs04} for more details). The sets
$S^\pm(X)$ restrict $d$ according to certain residue classes mod $Q$ in the case that
$Q$ is odd and squarefree, and $4Q$ in the case that $Q$ is even and squarefree.

\section{Moments of $L_E(1,\chi_d)$}
\label{section:moments}

Let
\begin{equation}
    M^\pm_E(X,k)=
    \frac{1}
    {\left| S^\pm(X) \right|}
    \sum_{d \in S^\pm(X)} L_E(1,\chi_d)^k.
\end{equation}
be the $k$th moment of $L_E(1,\chi_d)$.

The conjecture of Conrey-Farmer-Keating-Rubinstein-Snaith \cite[4.4]{kn:cfkrs}
says here that, for $k \geq 1$, $k \in \mathbb{Z}$,
\begin{equation}
    \label{eq:moment int}
    M^{\pm}_E(X,k)
    =
    \frac{1}{X}
    \int_0^X
    \Upsilon^{\pm}_k
    \left(\log(t)
    \right) dt
    +O(X^{-\tfrac{1}{2}+\epsilon})
\end{equation}
as $X\rightarrow\infty$, where $\Upsilon_k$ is the polynomial of degree
$k(k-1)/2$ given by the $k$-fold residue
\begin{eqnarray}
    \label{eq:upsilon}
    \Upsilon^{\pm}_k(x)&=&\frac{(-1)^{k(k-1)/2}2^{k}}{k!} \frac{1}{(2\pi
    i)^k}\\
&&\times    \oint \cdots \oint
    \frac{F^{\pm}_k(z_1,\ldots,z_k)\Delta(z_1^2,\ldots,z_k^2)^2}
    {\prod_{j=1}^k z_j^{2k-1}} e^{x\sum_{j=1}^kz_j}dz_1\ldots
    dz_k,\nonumber
\end{eqnarray}
where the contours above enclose the poles at $z_j=0$ and
\begin{equation}
    F^{\pm}_k(z_1,\ldots,z_k)=A^{\pm}_k(z_1,\ldots,z_k) \prod_{j=1}^k
    \left(
        \frac{\Gamma(1+z_j)}{\Gamma(1-z_j)}
        \left(\frac{Q}{4\pi^2} \right)^{z_j}
    \right)^{\tfrac{1}{2}}
    \prod_{1\leq i<j\leq k} \zeta(1+z_i+z_j).
\end{equation}
$A^{\pm}_k$, which depends on $E$, is the Euler product which is absolutely convergent for
$\sum_{j=1}^k |z_j|<1/2$,
\begin{equation}
    \label{eq:euler product A}
    A^{\pm}_k(z_1,\dots,z_k) =
    \prod_p F^\pm_{k,p}(z_1,\ldots,z_k)
    \prod_{1\le i < j \le k}
    \left(1-\frac{1}{p^{1+z_i+z_j}}\right)
\end{equation}
with, for $p \nmid Q$,
\begin{equation}
    \label{eq:Fkp big}
    F^\pm_{k,p}  =
    \left(1+\frac 1 p\right)^{-1}\left(\frac 1 p +\frac{1}{2}
    \left(\prod_{j=1}^k
    {\mathcal{L}}_p \left(\frac{1}{p^{1 +z_j}} \right)+
    \prod_{j=1}^k{\mathcal{L}}_p
    \left(\frac{-1}{p^{1 +z_j}}\right) \right)\right) .
\end{equation}
and, for $p\mid Q$,
\begin{equation}
    F^\pm_{k,p}=
    \prod_{j=1}^k
    {\mathcal{L}}_p \left(\frac{\pm a_p}{p^{1 +z_j}} \right).
\end{equation}
Because we are limiting ourselves to $Q$ squarefree
($Q$ prime in the $S^+$ case), we have $a_p=\pm 1$
when $p \mid Q$ and so
\begin{equation}
    \label{eq:Fkp}
    F^\pm_{k,p}=
    \begin{cases}
        \prod_{j=1}^k (1+p^{-1 -z_j})^{-1}  \quad \text{in the $S^-$ case, for $p \mid Q$}\\
        \prod_{j=1}^k (1-p^{-1 -z_j})^{-1}  \quad \text{in the $S^+$ case, for $p = Q$}.
    \end{cases}
\end{equation}

The r.h.s. of (\ref{eq:moment int}) is~\cite{kn:cfkrs} asymptotically, as $X\rightarrow \infty$,
\begin{equation}
    \label{eq:leading asymp}
    M^{\pm}_E(X,k)\sim A^{\pm}(k)M_O(\lfloor \log X \rfloor, k)
\end{equation}
where
\begin{eqnarray}
    \label{eq:Ak}
    &&A^{\pm}(k)= \\&&\qquad\prod_{p \nmid Q} \left( 1-p^{-1}\right)^{k(k-1)/2}
          \left(\frac{p}{p+1}\right)
          \left(
              \frac1p + \frac12
              \left(
                  {\mathcal{L}}_p(1/p)^{k} +
                  {\mathcal{L}}_p(-1/p)^{k}
              \right)
          \right) \nonumber\\ \notag
        &&\qquad  \times\prod_{p \mid Q} \left( 1-p^{-1}\right)^{k(k-1)/2}
          {\mathcal{L}}_p(\pm a_p/p)^{k}
\end{eqnarray}
with
\begin{equation}
    \label{eq:M_O}
    M_O(N,k)=2^{2Nk}\prod_{j=1}^N\frac{{\rm
    \Gamma}(N+j-1){\rm \Gamma}(k+j-1/2)}{{\rm \Gamma}(j-1/2){\rm
    \Gamma}(k+j+N-1)}.
\end{equation}
The leading asymptotics given above for the moments of $L_E(1,\chi_d)$
was first made in~\cite{kn:keasna00b} and~\cite{kn:confar00}, though the
arithmetic factor was off for primes dividing $Q$. One nice thing
about~(\ref{eq:leading asymp}) is that it makes sense for complex
values of $k$ and in~\cite{kn:keasna00b} was conjectured to hold for
$\Re{k}> -1/2$.

In \cite{kn:ckrs00} it is shown how the conjectured asymptotics
for moments can be used to obtain information concerning the distribution
of values of $L_E(1,\chi_d)$. That paper discusses the importance of the first pole of
the r.h.s. of~(\ref{eq:M_O}) at $k=-1/2$ in analyzing the number of vanishings
of $L_E(1,\chi_d)$.

\section{Vanishings of $L_E(1,\chi_d)$ in progressions}
\label{section:vanishings}

We fix a prime $q \nmid Q$ and restrict $d$ further according to residue classes
mod $q$ as follows. For $\lambda=\pm 1$ we set
\begin{equation}
     S^\pm(X;q,\lambda) = \{ d \in S^\pm(X); \chi_d(q) =\lambda \}
\end{equation}

Let
\begin{equation}
    R^{\pm}_q(X)=
    \frac{
        \sum_{d \in S^{\pm}(X;q,1) \atop L_E(1,\chi_d)=0}  1
    }
    {
        \sum_{d \in S^{\pm}(X;q,-1) \atop L_E(1,\chi_d)=0}  1
    }
\end{equation}
denote the number of ratio of the number of vanishings of $L_E(1,\chi_d)$, with
$d \in S^{\pm}$, sorted according to residue classes mod $q$.

To study this ratio we need to look at the
moments:
\begin{equation}
    \label{eq:moments subfamilies}
    M^\pm_E(X,k;q,\lambda)=
    \frac{1}
    {\left| S^\pm(X;q,\lambda) \right|}
    \sum_{d \in S^\pm(X;q,\lambda)} L_E(1,\chi_d)^k.
\end{equation}
The conjecture in~\cite{kn:cfkrs} then gives
\begin{equation}
    \label{eq:conjecture subfamilies}
    M^\pm_E(X,k;q,\lambda)
    =
    \frac{1}{X}
    \int_0^X
    \Upsilon^{\pm}_{k,q,\lambda}
    \left(\log(t)
    \right) dt
    +O(X^{-\tfrac{1}{2}+\epsilon})
\end{equation}
where $\Upsilon^{\pm}_{k,q,\lambda}(x)$ is given by the same formula as
in (\ref{eq:upsilon}) but with a slight but important modification: the local factor
corresponding to the prime $q$, $F^\pm_{k,q}$, gets replaced by
\begin{equation}
    \label{eq:F lambda}
    F^\pm_{k,q,\lambda} =
    \prod_{j=1}^k (1-\lambda a_q q^{-1 -z_j}+q^{-1-2 z_j})^{-1}.
\end{equation}
Similarly, in (\ref{eq:Ak}), the local factor
\begin{equation}
    \left(\frac{q}{q+1}\right)
    \left(
        \frac1q + \frac12
        \left(
            { \mathcal{L}}_q(1/q)^{k} +
             {\mathcal{L}}_q(-1/q)^{k}
        \right)
    \right)
\end{equation}
at the prime $q$ gets replaced by
\begin{equation}
    {L}_q(\lambda/q)^{k} = (1-\lambda a_q q^{-1}+q^{-1})^{-k}.
\end{equation}
From this we immediately surmise several things. First, $R^\pm_q(X)$
which is conjectured to be, asymptotically, equal to the ratio of the residues of
the two moments~(\ref{eq:conjecture subfamilies}), corresponding to
$\lambda=1$ and $-1$, at the pole $k=-1/2$ should thus equal, up to
leading order,
\begin{equation}
    \left( \frac{q+1-a_q}{q+1+a_q} \right)^{1/2}.
\end{equation}
Second, when $a_q=0$, the complete asymptotic expansion for both
moments are identical up to the conjectured error of size $O(X^{-1/2+\epsilon})$.
The reason for this is that, in~(\ref{eq:F lambda}), if $a_q=0$, there is no
dependence on $\lambda$. Indulging in conjectural bravado, we predict
that when $a_q=0$
\begin{equation}
    R^\pm_q(X) = 1 + O (X^{-1/2+\epsilon})
\end{equation}
and similarly for $R_q(X)$ in (\ref{eq:Rq}). This fits well with our numeric data.
See section~\ref{section:numerics} and also Table~1 in \cite{kn:ckrs00} .

Third, from this formula for the moments we are able to work out, in principle,
arbitrarily many terms in the asymptotic expansion of $R^\pm_q(X)$. Below, we describe
the next to leading term in detail. It is of size $O(1/\log(X))$.
The lower terms in the asymptotics of $R^{\pm}_q(X)$
do depend on whether we are looking at $S^{+}(X)$ as opposed
to $S^{-}(X)$. This arises from the fact that the
local factors $F^\pm_{k,p}$ for $p \mid Q$ in equation~(\ref{eq:Fkp})
depend on whether we are looking at $S^+$ or $S^-$.
While this does not affect the main term $R_q$, it does show up in the secondary terms.

\section{Evaluating the first two terms of $M^\pm_E(X,k;q,\lambda)$}
\label{section:two terms of moments}

To evaluate the residue that defines $\Upsilon^{\pm}_{k,q,\lambda}$ we
need to examine the multiple Laurent series about $z_j=0$ of the corresponding
integrand.
In the numerator, we must evaluate the coefficient of
$\prod_{j=1}^k z_j^{2k-2}$ of degree $2k(k-1)$. Now
$\Delta(z_1^2,\ldots,z_k^2)^2$ is a homogeneous polynomial consisting of terms of
degree $4 {k \choose 2} = 2k(k-1)$. However, the poles of
$\prod_{1\leq i<j\leq k} \zeta(1+z_i+z_j)$ cancel ${k \choose 2}$
factors of the Vandermonde. Therefore, in computing the residue, we
only need to take terms from the series for $e^{x\sum_{j=1}^kz_j}$ up to
degree ${k \choose 2}$. From this we see that $\Upsilon^{\pm}_{k,q,\lambda}(x)$
is a polynomial in $x$ of degree ${k \choose 2}$.

To obtain the leading two terms of $\Upsilon^{\pm}_{k,q,\lambda}(x)$, i.e.
those of degree ${k \choose 2}$ and ${k \choose 2}-1$ in $x$, we need to evaluate
the constant and linear terms in the multiple Maclaurin series
of the function
\begin{eqnarray}
    h^\pm_k(z;q,\lambda) &=&
    A^{\pm}_k(z_1,\ldots,z_k;q,\lambda) \prod_{j=1}^k
    \left(
        \frac{\Gamma(1+z_j)}{\Gamma(1-z_j)}
        \left(\frac{Q}{4\pi^2} \right)^{z_j}
    \right)^{\tfrac{1}{2}}\\
    &&\qquad\qquad \nonumber\times\prod_{1\leq i<j\leq k} \zeta(1+z_i+z_j) (z_i+z_j).
\end{eqnarray}
Here $A^{\pm}_k(z_1,\ldots,z_k;q,\lambda)$ is the same as the function
$A^{\pm}_k(z_1,\ldots,z_k)$ but with the local factor $F^\pm_{k,q}$ replaced by
$F^\pm_{k,q,\lambda}$.

For example, the term involving $x^{k(k-1)/2}$ of $\Upsilon^{\pm}_{k,q,\lambda}(x)$
is equal to
\begin{eqnarray}
   && h^\pm_k(0;q,\lambda)
    \frac{(-1)^{k(k-1)/2}2^{k}}{k!} \frac{1}{(2\pi i)^k}\\
   &&\qquad\qquad\times \oint \cdots \oint
    \frac{\Delta(z_1^2,\ldots,z_k^2)^2}
    {\prod_{j=1}^k z_j^{2k-1}}
    \frac{e^{x\sum_{j=1}^kz_j}}
    {\prod_{1\leq i<j\leq k} (z_i+z_j)}
    dz_1\ldots dz_k.\nonumber
\end{eqnarray}
It is shown in~\cite{kn:cfkrs1} that the above equals
\begin{equation}
    \label{eq:leading term}
    h^\pm_k(0;q,\lambda)
    g_k(O^+)
    x^{k(k-1)/2}
\end{equation}
where
\begin{equation}
   \label{eq:g_k}
   g_k(O^+)=2^{k(k+1)/2}\prod_{j=1}^{k-1}\frac{j!}{2j!}.
\end{equation}
We also have
\begin{equation}
    \label{eq:h 0}
    h^\pm_k(0;q,\lambda)
    = A^{\pm}_k(0,\ldots,0;q,\lambda).
\end{equation}

To compute the leading two terms of the moments we prefer to write
\begin{equation}
    h^\pm_k(z;q,\lambda) = \exp(\log{h^\pm_k(z;q,\lambda)})
\end{equation}
and evaluate the constant and linear terms of
\begin{equation}
    \label{eq:h series}
    \log{h^\pm_k(z;q,\lambda)}=\alpha_k^\pm(q,\lambda) + \beta^\pm_k(q,\lambda) \sum z_j + \ldots.
\end{equation}
Notice that the linear terms all share the same coefficient because $h^\pm_k(z;q,\lambda)$
is symmetric in the $z_j$'s.

The constant term can be pulled out of the integral as $e^{\alpha_k^\pm(q,\lambda)}=h^\pm_k(0;q,\lambda)$.
The linear terms can be absorbed into the $\exp(x\sum_{j=1}^k z_j)$.
Dropping the terms of degree two or higher in $\log{h^\pm_k(z;q,\lambda)}$ we
can evaluate the residue using~(\ref{eq:leading term}):
\begin{equation}
    \label{eq:leading term 2}
    h^\pm_k(0;q,\lambda)
    g_k(O^+)
    (x+\beta^\pm_k(q,\lambda))^{k(k-1)/2}
\end{equation}
and thus find that
\begin{equation}
    \label{eq:upsilon first 2 terms}
    \Upsilon^{\pm}_{k,q,\lambda}(x)
    = h^\pm_k(0;q,\lambda)
    g_k(O^+)
    (x^{\frac{k(k-1)}{2}} + \frac{k(k-1)}{2} \beta^\pm_k(q,\lambda) x^{\frac{k(k-1)}{2}-1} + \ldots).
\end{equation}
Inserting (\ref{eq:upsilon first 2 terms}) into (\ref{eq:conjecture subfamilies})
and integrating, we obtain
\begin{eqnarray}
   M^\pm_E(X,k;q,\lambda)
    &=&
    \frac{ h^\pm_k(0;q,\lambda) g_k(O^+)} {X}\\
  &&\times  \int_0^X
    \left(
         \log(t)^{\frac{k(k-1)}{2}} +
         \frac{k(k-1)\beta^\pm_k(q,\lambda)}{2}
         \log(t)^{\frac{k(k-1)}{2}-1}
    \right) dt \nonumber\\
   &&\qquad\qquad +O(\log(X)^{\frac{k(k-1)}{2}-2}) \notag
\end{eqnarray}
and hence
\begin{eqnarray}
    \label{eq:conjecture moments 2 terms}
    &&M^\pm_E(X,k;q,\lambda)
    =
    h^\pm_k(0;q,\lambda) g_k(O^+) \log(X)^{\frac{k(k-1)}{2}}\\
  &&\qquad\qquad \times  \left(
        1+\frac{k(k-1)}{2 \log(X)}(\beta^\pm_k(q,\lambda)-1)
    \right)+O(\log(X)^{\frac{k(k-1)}{2}-2}).\nonumber
\end{eqnarray}

Therefore, the remaining work is to compute above the coefficient
$\beta^\pm_k(q,\lambda)$. To do so we evaluate individually the linear
terms in the Maclaurin expansions of:
\begin{equation}
    \label{eq:log gamma}
    \frac{1}{2}
    \log \prod_{j=1}^k
    \left(
        \frac{\Gamma(1+z_j)}{\Gamma(1-z_j)}
        \left(\frac{Q}{4\pi^2} \right)^{z_j}
    \right),
\end{equation}
\begin{equation}
    \label{eq:log zeta}
    \log \prod_{1\leq i<j\leq k} \zeta(1+z_i+z_j) (z_i+z_j),
\end{equation}
and
\begin{equation}
    \label{eq:log A}
    \log A^{\pm}_k(z_1,\ldots,z_k;q,\lambda).
\end{equation}

First, $\log\Gamma(1+z) = -\gamma z +\frac{\pi^2}{12}z^2 + \ldots$ hence
\begin{equation}
    \frac{1}{2}
    \log
    \left(
        \frac{\Gamma(1+z)}{\Gamma(1-z)}
        \left(\frac{Q}{4\pi^2} \right)^{z}
    \right) = (-\gamma + \log(Q^{1/2}/(2\pi))) z + \ldots
\end{equation}
and so (\ref{eq:log gamma}) equals
\begin{equation}
    (-\gamma + \log(Q^{1/2}/(2\pi))) \sum z_j + \ldots.
\end{equation}
Next,
\begin{equation}
    \zeta(1+z_i+z_j) (z_i+z_j) = 1 + \gamma(z_i+z_j) + \ldots
\end{equation}
so
\begin{eqnarray}
    \prod_{1\leq i<j\leq k} \zeta(1+z_i+z_j) (z_i+z_j)
    &=& 1 + \gamma \sum_{1\leq i<j\leq k} (z_i+z_j) +
    \ldots\nonumber \\
    &= &1 + (k-1)\gamma \sum z_j + \ldots
\end{eqnarray}
Therefore, (\ref{eq:log zeta}) equals
\begin{equation}
    (k-1) \gamma \sum z_j + \ldots.
\end{equation}
We now turn to (\ref{eq:log A}). The function
$A^{\pm}_k(z_1,\ldots,z_k;q,\lambda)$ is given by~(\ref{eq:euler product A}) except
that the local factor
at $p=q$, namely $F^\pm_{k,q}$, gets replaced by~(\ref{eq:F lambda}).
To find the coefficient of $\sum z_j$ in the Maclaurin series for
\begin{equation}
    \prod_{1\le i < j \le k}
    \left(1-\frac{1}{p^{1+z_i+z_j}}\right)
\end{equation}
we can, because the above is symmetric in the $z_j$'s,
differentiate with respect to $z_1$ and set all $z_j$ equal to 0.
We thus find that the coefficient of $\sum z_j$ equals
\begin{equation}
    \frac{(k-1) \log p}{p-1}.
\end{equation}

Next we consider the contribution from the local factor when $p=q$:
\begin{equation}
    \log F^\pm_{k,q,\lambda} =
    - \sum_{j=1}^k \log(1-\lambda a_q q^{-1 -z_j}+q^{-1-2 z_j}).
\end{equation}
Differentiating w.r.t.~$z_1$ and setting all $z_j=0$ we find that the
coefficient of $\sum z_j$ in the Maclaurin series for
$\log F^\pm_{k,q,\lambda}$ equals
\begin{equation}
    \frac{\log q (\lambda a_q-2)}{\lambda a_q - q -1}.
\end{equation}

Finally, we consider the local factor when $p\neq q$. If
$p \mid Q$, we have, on taking the logarithm of (\ref{eq:Fkp}),
differentiating w.r.t.~$z_1$, setting all $z_j=0$, that the
coefficient of $\sum z_j$ in the series for $\log F^\pm_{k,p}$ equals
\begin{equation}
    \begin{cases}
          \log(p)/(1+p) \quad \text{in the $S^-$ case}\\
          \log(p)/(1-p) \quad \text{in the $S^+$ case}.
    \end{cases}
\end{equation}
If $p\nmid Q$, taking the logarithm of (\ref{eq:Fkp big}),
differentiating w.r.t.~$z_1$, and letting $z_j=0$, we get the coefficient
of $\sum z_j$ equal to
\begin{equation}
    \log(p)
    \left(
        \frac{(2-a_p)f_1(p)^{-k-1} + (2+a_p)f_2(p)^{-k-1}}
             {2+p\left(f_1(p)^{-k} + f_2(p)^{-k} \right)}
    \right)
\end{equation}
where
\begin{eqnarray}
    \label{eq:f1 f2}
    f_1(p)=1-a_p/p+1/p \notag \\
    f_2(p)=1+a_p/p+1/p.
\end{eqnarray}

Hence, adding all the coefficients of $\sum z_j$ we find that
$\beta^\pm_k(q,\lambda)$ in (\ref{eq:h series}), and hence in
(\ref{eq:upsilon first 2 terms}), equals
\begin{equation}
    \label{eq:beta explicit}
    (k-2) \gamma +
    \log(Q^{1/2}/(2\pi))
    + \sum_{p} \beta_k(p)
\end{equation}
where
\begin{equation}
    \beta_k(p) = \frac{(k-1) \log p}{p-1} +
    \begin{cases}
        \frac{\log(q) (\lambda a_q-2)}{\lambda a_q - q -1}
        \quad \text{if $p=q$}\\
        \log(p)
        \left(
            \frac{(2-a_p)f_1(p)^{-k-1} + (2+a_p)f_2(p)^{-k-1}}
                 {2+p\left(f_1(p)^{-k} + f_2(p)^{-k} \right)}
        \right) \quad \text{if $p \neq q$, $p \nmid Q$} \\
        \log(p)/(1+p) \quad \text{if $p \mid Q$, in the $S^-$ case}\\
        \log(p)/(1-p) \quad \text{if $p \mid Q$, in the $S^+$ case.}
    \end{cases}
\end{equation}
Notice that the only dependence in $\beta^\pm_k(q,\lambda)$ on $q$
is in the term
\begin{equation}
    \label{eq:q dependency}
    \beta_k(q) = \frac{(k-1) \log q}{q-1} +
    \frac{\log(q) (\lambda a_q-2)}{\lambda a_q - q -1}.
\end{equation}

\section{Conjecture for the first two terms in $R^\pm_q(X)$}
\label{section:conjecture}

Dividing $M^\pm_E(X,k;q,1)$ by $M^\pm_E(X,k;q,-1)$, using
equation~(\ref{eq:conjecture moments 2 terms})
\begin{equation}
    \frac{M^\pm_E(X,k;q,1)}{M^\pm_E(X,k;q,-1)}
    =
    \frac{h^\pm_k(0;q,1)}{h^\pm_k(0;q,-1)}
    \frac{
        \left(
            1+\frac{k(k-1)}{2 \log(X)}(\beta^\pm_k(q,1)-1)
        \right)
    }
    {
        \left(
            1+\frac{k(k-1)}{2 \log(X)}(\beta^\pm_k(q,-1)-1)
        \right)
    }
    +O(\log(X)^{-2}).
\end{equation}
The first factor $\frac{h^\pm_k(0;q,1)}{h^\pm_k(0;q,-1)}$ equals
\begin{equation}
    \left( \frac{q+1-a_q}{q+1+a_q} \right)^{-k}.
\end{equation}
Interpolating to $k=-1/2$ gives our conjecture:
\begin{conjecturemr}
For $q\nmid Q$
\begin{equation}
    \label{eq:second term}
    R^{\pm}_q(X) =
    R_q
    \frac{1+\frac{3}{8 \log(X)}(\beta^\pm_{-\frac{1}{2}}(q,1)-1)}
    {1+\frac{3}{8 \log(X)}(\beta^\pm_{-\frac{1}{2}}(q,-1)-1)}
    +O(\log(X)^{-2})
\end{equation}
where $\beta^\pm_{-\frac{1}{2}}(q,\lambda)$ is given explicitly by
equation~(\ref{eq:beta explicit}).
The implied constant in the remainder term
depends on $E$ and $q$, and thus also on $a_q$.
\end{conjecturemr}

\section{Numerical Data}
\label{section:numerics}

We verify the conjecture described above for over two thousand elliptic curves
and the sets $S^{\pm}_E(X)$, with $X=10^8$. Altogether we have 2398 datasets.
The curves in question and the method for computing $L_E(1,\chi_d)$ are
detailed in~\cite{kn:ckrs04}. Tables of $L$-values can be obtained
from ~\cite{kn:rub04}.

\begin{figure}[htp]
    \centerline{
            \psfig{figure=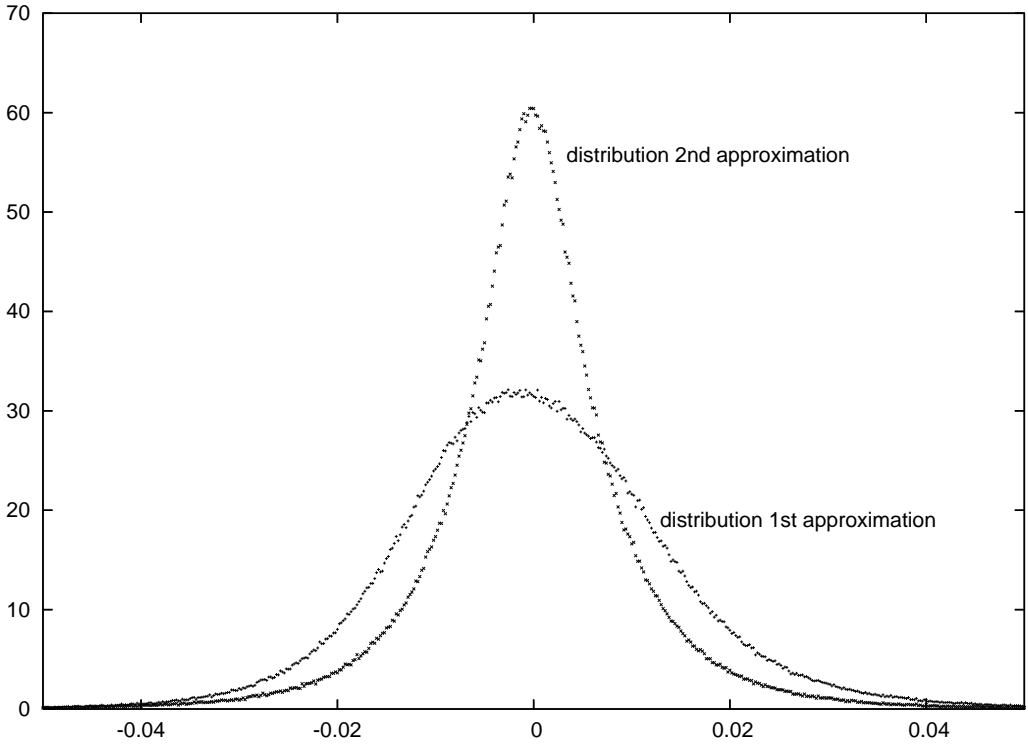,width=4in,angle=-90}
    }
    \caption
    {Distribution first approximation v.s. second approximation for ratio of vanishings }
    \label{fig:dist 1st vs 2nd}
\end{figure}

We first depict in Figure~\ref{fig:dist 1st vs 2nd}
the distribution of the remainder in comparing $R^{\pm}_q(X)$
to the conjectured first and second order approximations.
More precisely, for our 2398 datasets, we examine the
distribution of values of
\begin{equation}
    \label{eq:reality minus conj1}
    R^{\pm}_q(X)- R_q
\end{equation}
and of
\begin{equation}
    \label{eq:reality minus conj2}
    R^{\pm}_q(X)-R_q
    \frac{1+\frac{3}{8 \log(X)}(\beta^\pm_{-\frac{1}{2}}(q,1)-1)}
    {1+\frac{3}{8 \log(X)}(\beta^\pm_{-\frac{1}{2}}(q,-1)-1)}
\end{equation}
with $X=10^8$, $q \leq 3571$.
We break up the horizontal axis into small bins of size $.0002$ and
count how often the values fall within a given bin.
The difference in (\ref{eq:reality minus conj2}) has smaller variance reflecting an overall
better fit of the second order
approximation compared with the first. These distributions are not Gaussian.
There are yet further lower terms and these are given
by complicated sums involving the Dirichlet coefficients of $L_E(s)$, and $q$.

\begin{figure}[htp]
    \centerline{
            \psfig{figure=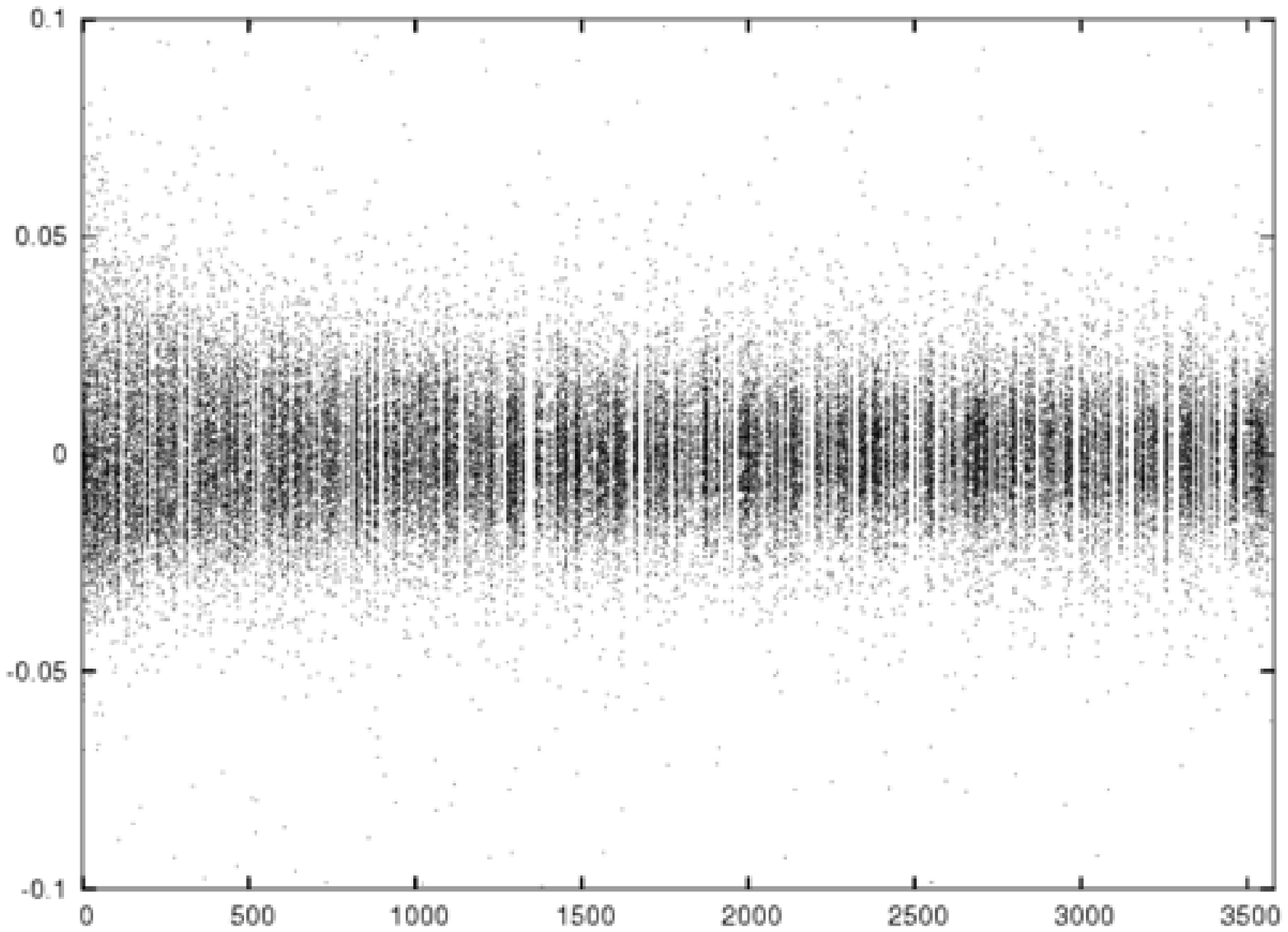,width=4in,angle=-90}
    }
    \centerline{
            \psfig{figure=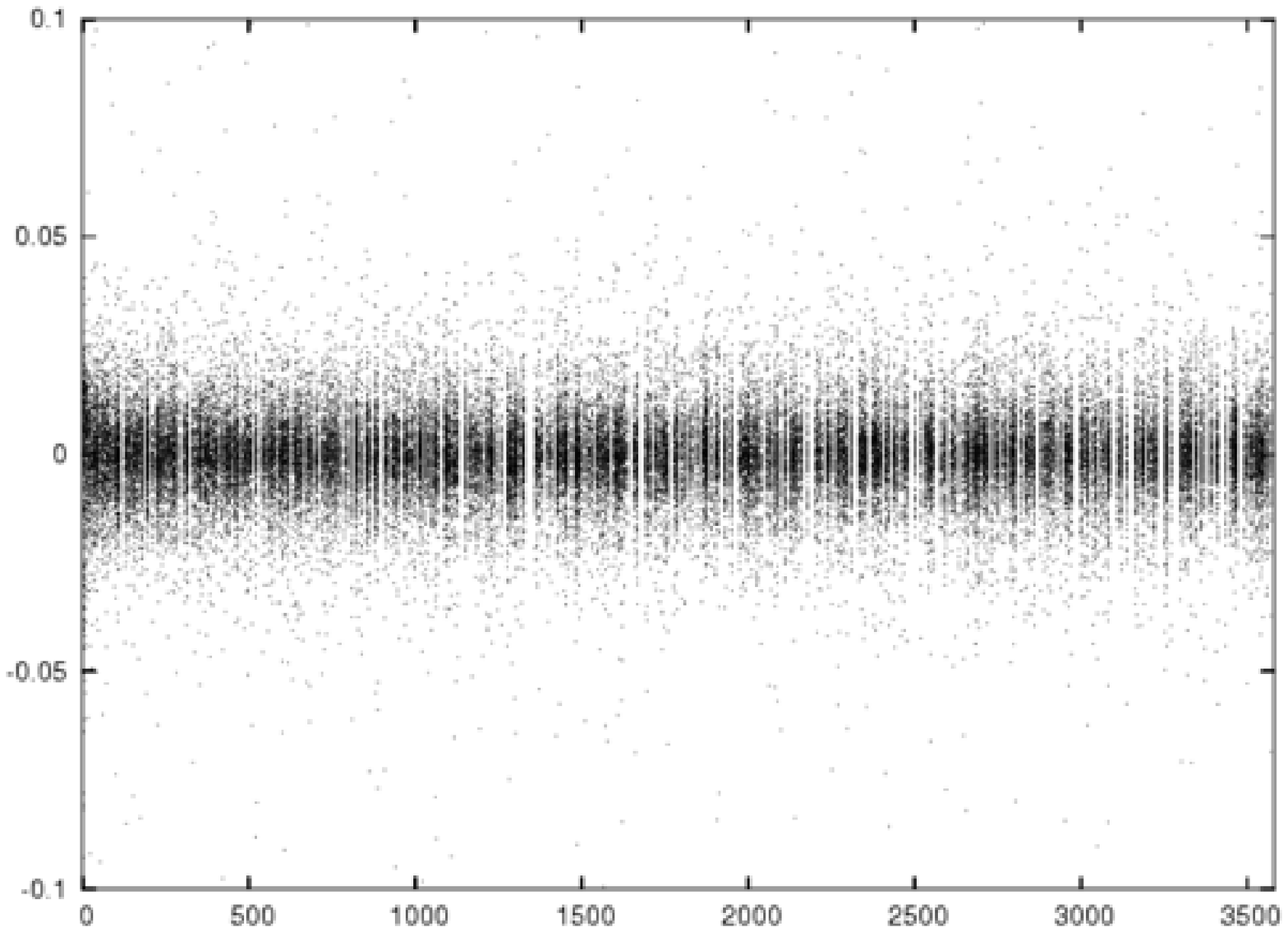,width=4in,angle=-90}
    }
    \caption
    {A plot for one hundred datasets of $R^{\pm}_q(10^8)-R_q$, top plot,
     and of (\ref{eq:reality minus conj2}), bottom plot, for $2\leq q\leq 3571$.}
    \label{fig:1st vs 2nd}
\end{figure}

In the first plot of Figure~\ref{fig:1st vs 2nd} we depict,
for one hundred of our datasets,
the raw data for the values given by equation~(\ref{eq:reality minus conj1}).
The horizontal axis is $q$.
For each $q$ on the horizontal axis there are 100 points corresponding to the
100 values, one for each dataset, of $R^{\pm}_q(X)-R_q$, with $X=10^8$.
We see the values fluctuating
about zero, most of the time agreeing to within about $.02$.
The convergence in $X$ is predicted from the secondary term to be
logarithmically slow and one gets a better fit by including the second order
term.

This is depicted in the second plot of Figure~\ref{fig:1st vs 2nd}
which shows the difference given in~(\ref{eq:reality minus conj2}).
again with $X=10^8$, and the same one hundred elliptic curves $E$.
We see an improvement to the first plot which uses just the main term.
We only depict data for 100 datasets in these plots since otherwise there would
be too many data points leading to a thick black mess.

\begin{figure}[htp]
    \centerline{
            \psfig{figure=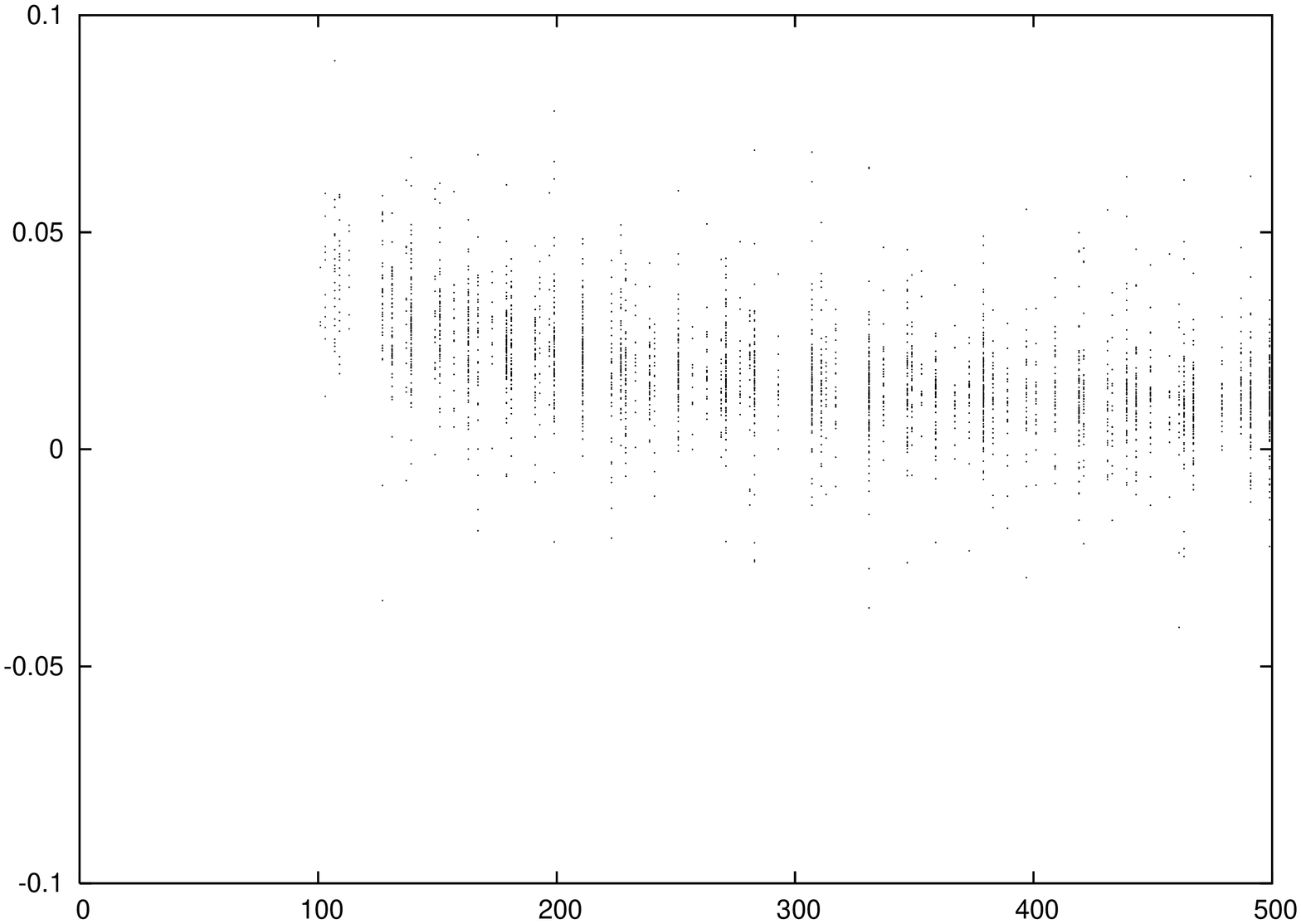,width=1.25in,angle=-90}
            \psfig{figure=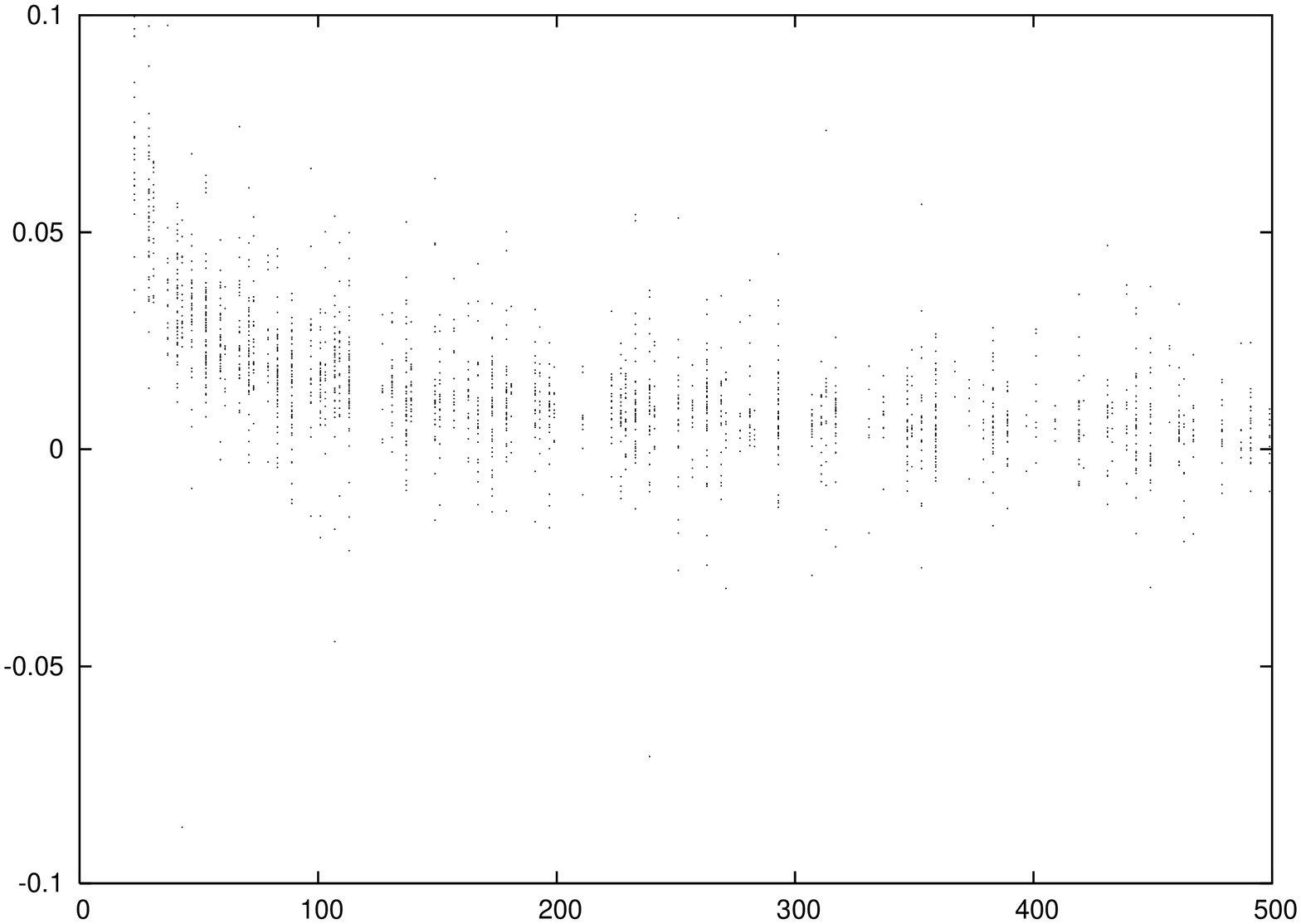,width=1.25in,angle=-90}
            \psfig{figure=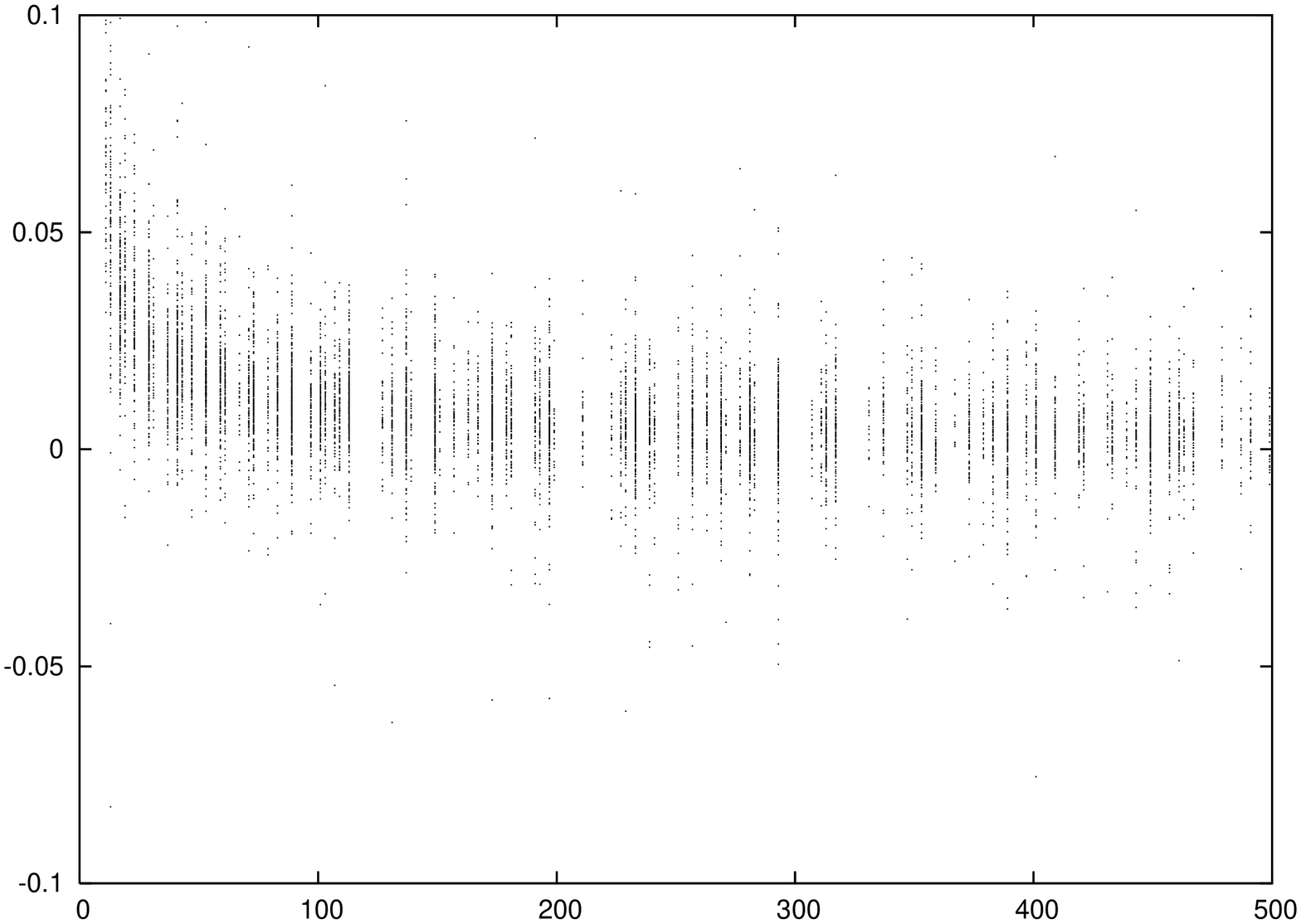,width=1.25in,angle=-90}
    }
    \centerline{
            \psfig{figure=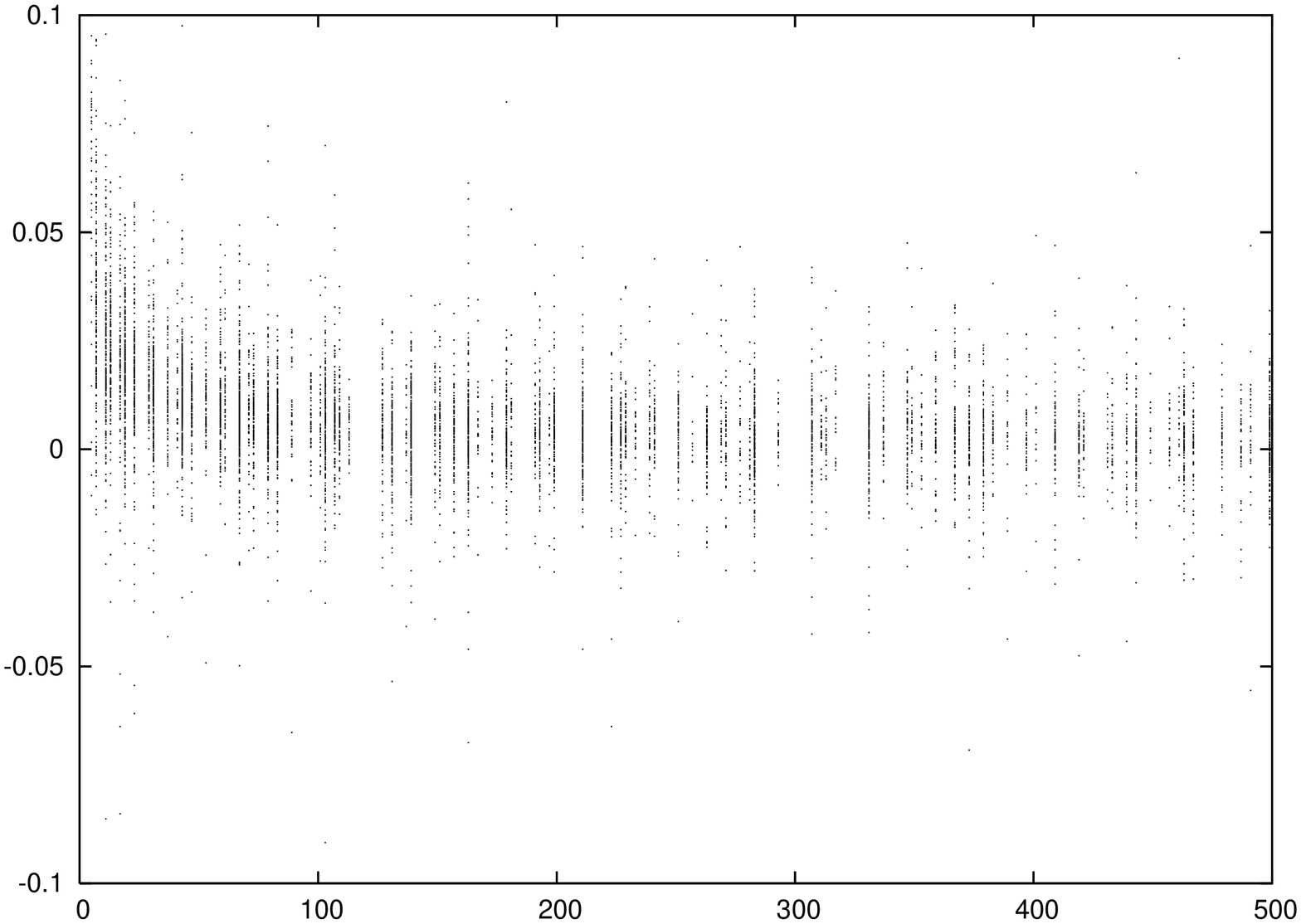,width=1.25in,angle=-90}
            \psfig{figure=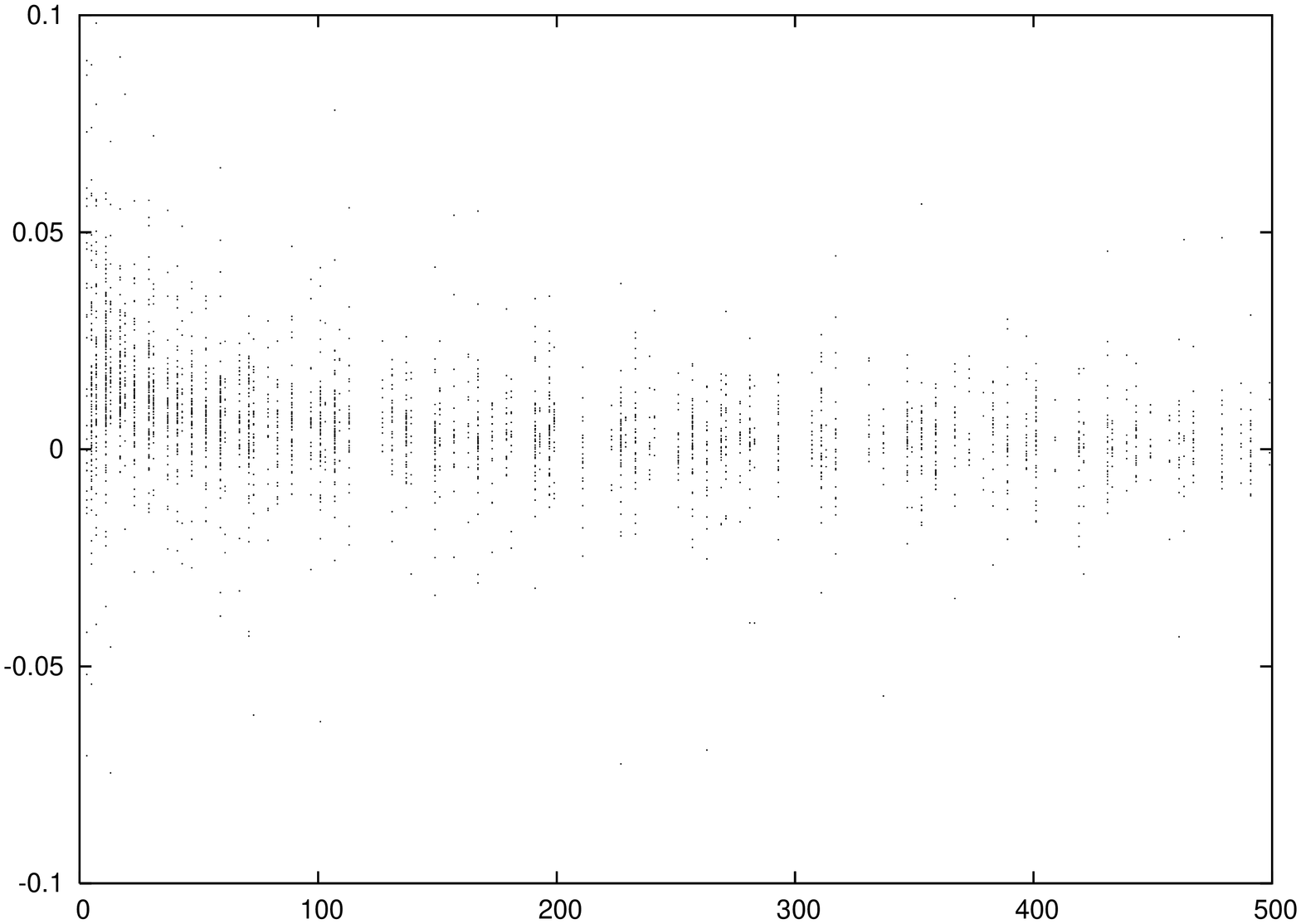,width=1.25in,angle=-90}
            \psfig{figure=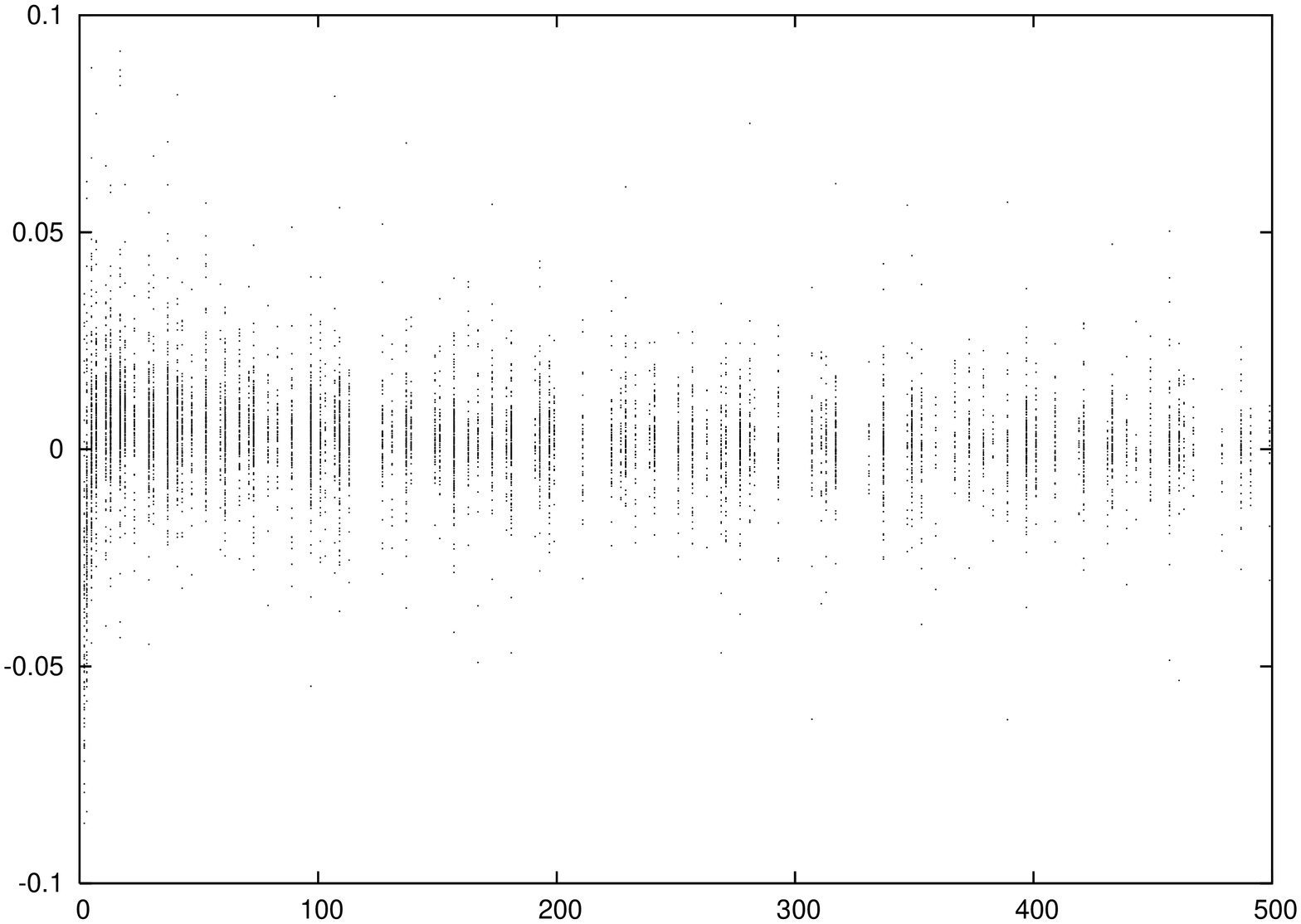,width=1.25in,angle=-90}
    }
    \centerline{
            \psfig{figure=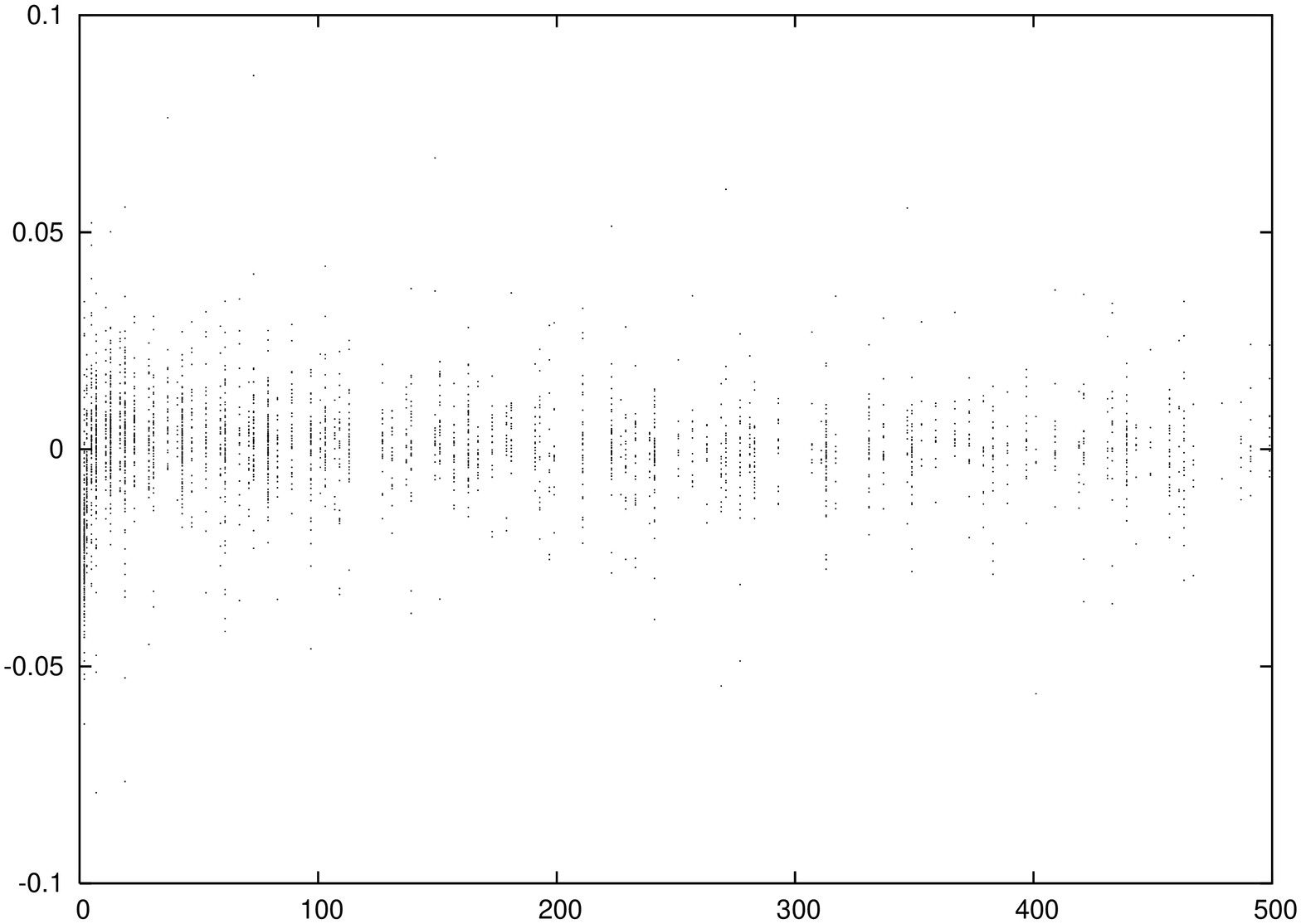,width=1.25in,angle=-90}
            \psfig{figure=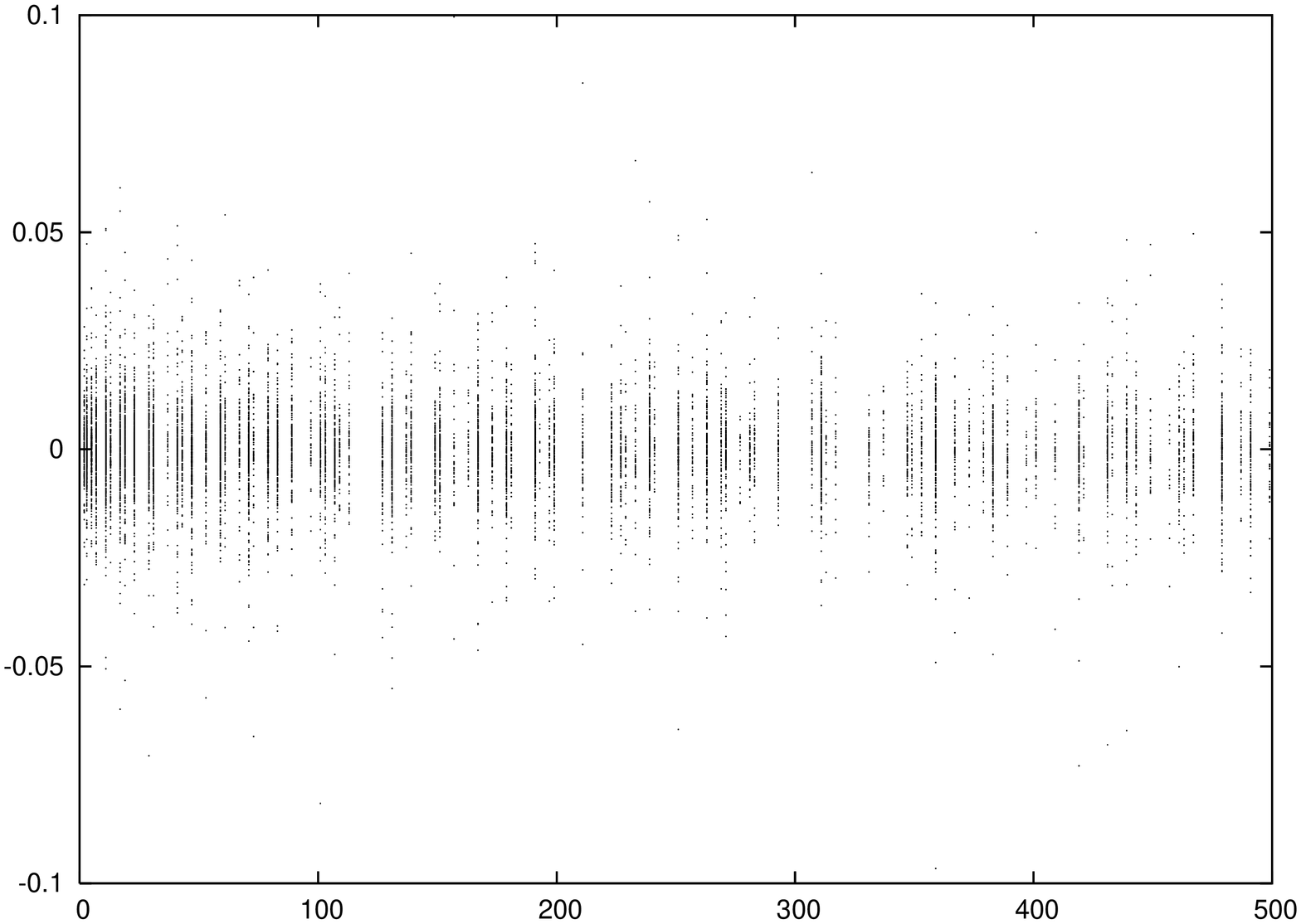,width=1.25in,angle=-90}
            \psfig{figure=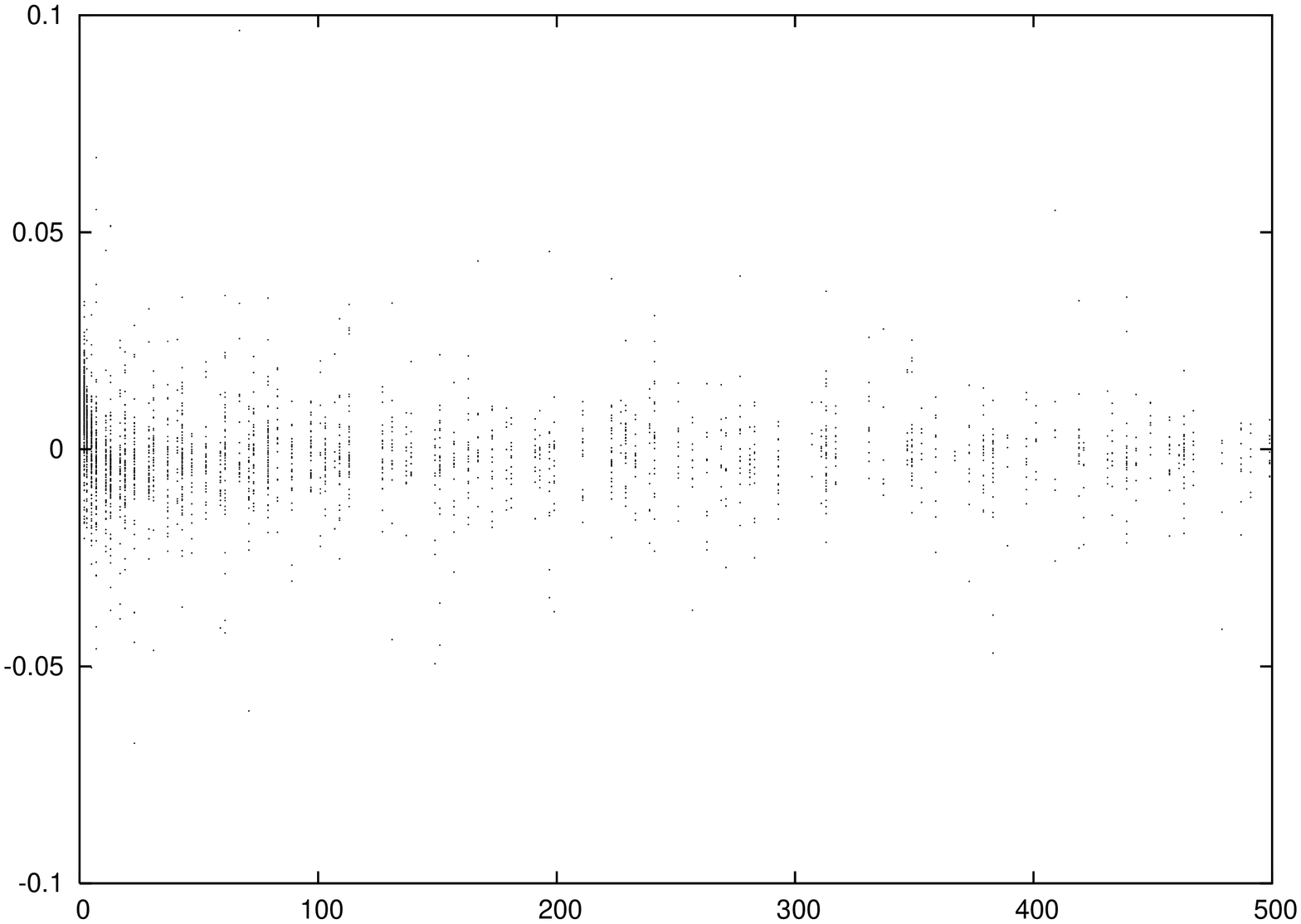,width=1.25in,angle=-90}
    }
    \centerline{
            \psfig{figure=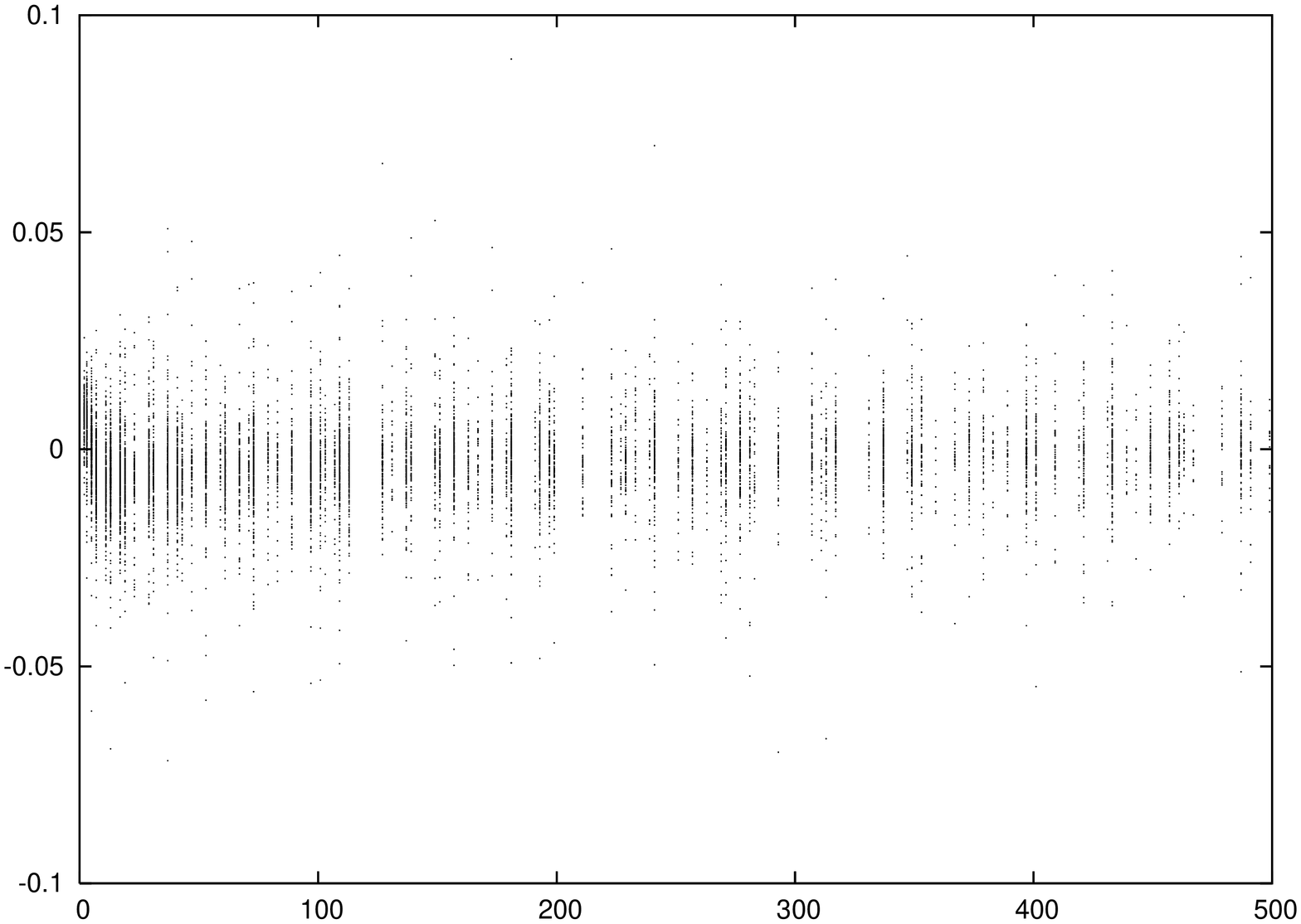,width=1.25in,angle=-90}
            \psfig{figure=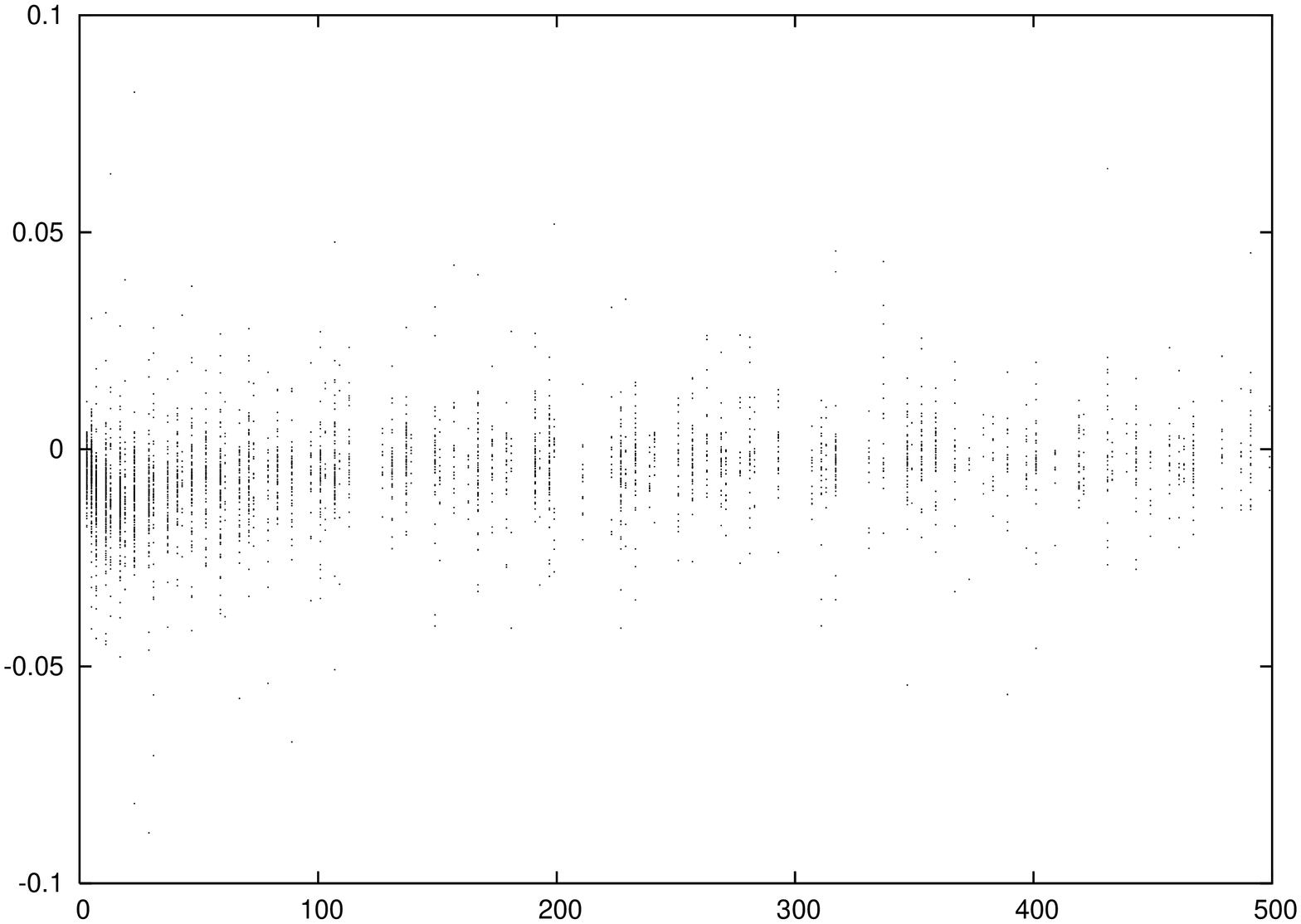,width=1.25in,angle=-90}
            \psfig{figure=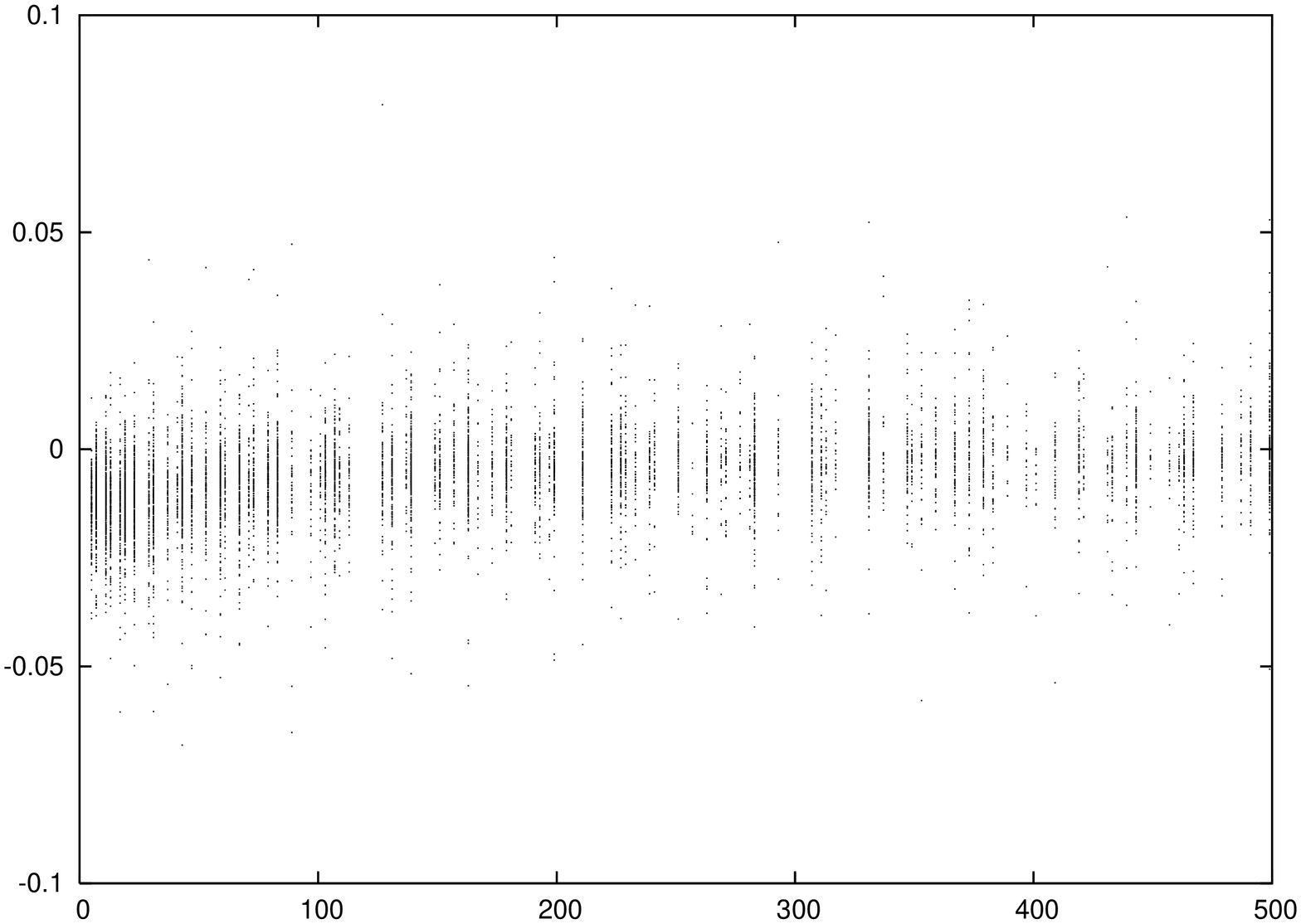,width=1.25in,angle=-90}
    }
    \centerline{
            \psfig{figure=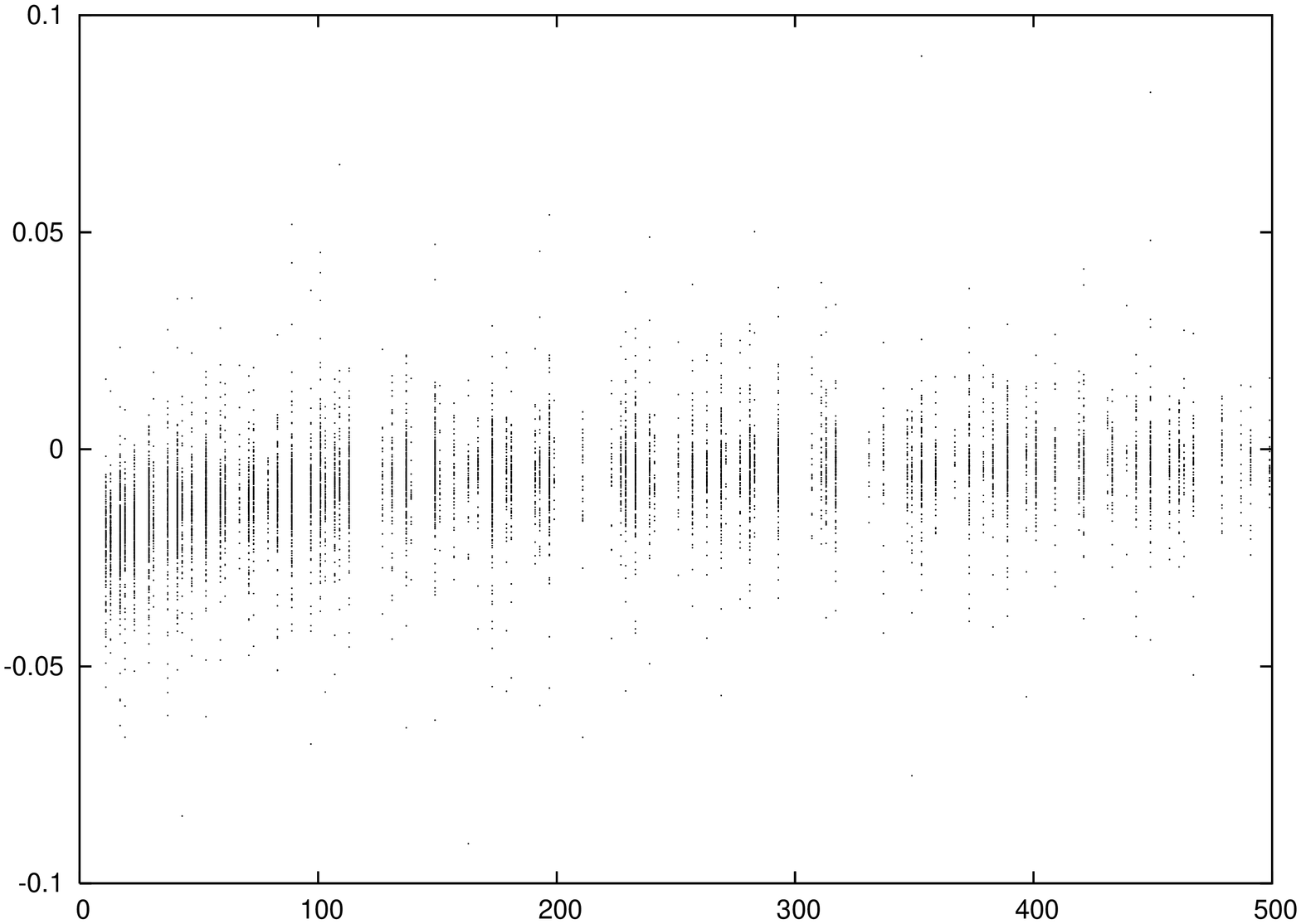,width=1.25in,angle=-90}
            \psfig{figure=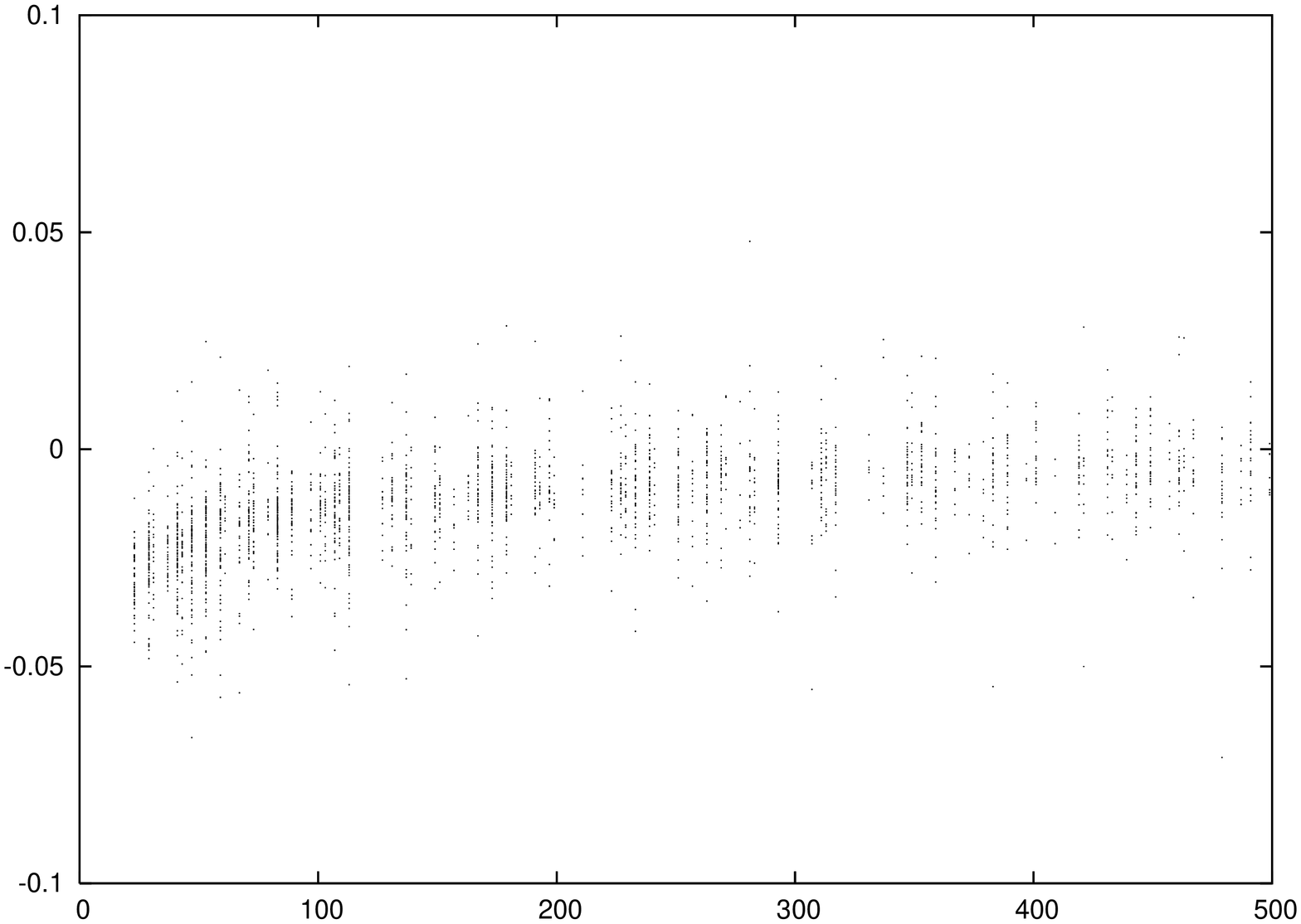,width=1.25in,angle=-90}
            \psfig{figure=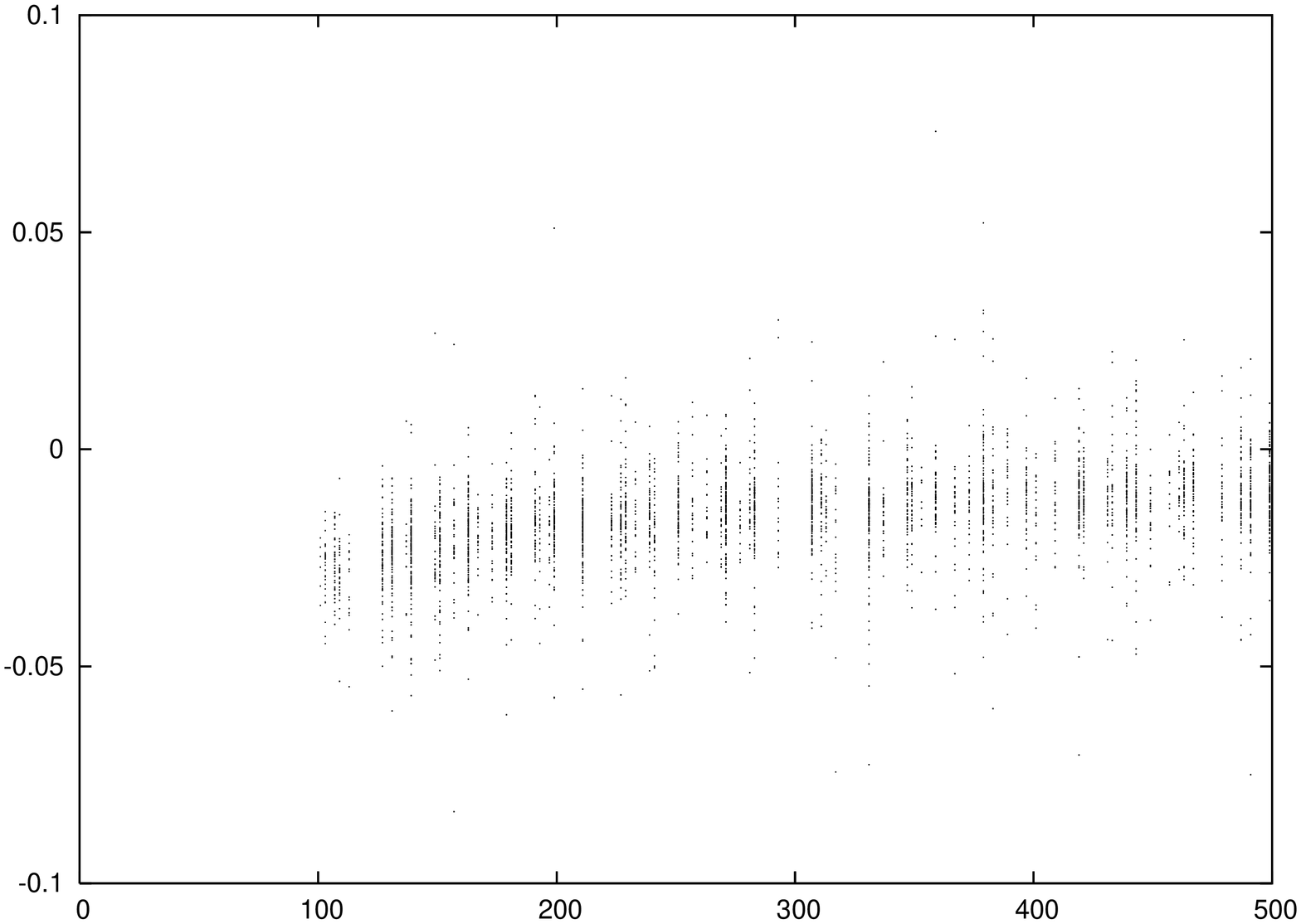,width=1.25in,angle=-90}
    }
    \caption
    {Left to right, top to bottom: $n=-20,-9,-6,-4,-3,-2,-1,0,1,2,3,4,6,9,20$.
     Values of $R^{\pm}_q(X)-R_q$, with $X=10^8$, $2\leq q <500$, for the
     subset of our elliptic curves satisfying $a_q=n$.
     The blank white area on the left of the plots for larger $n$ reflects Hasse's theorem
     that $|a_q| < 2 q^{1/2}$ which restricts how small $q$ can be given $a_q=n$.}
    \label{fig:seq 1st}
\end{figure}

\begin{figure}[htp]
    \centerline{
            \psfig{figure=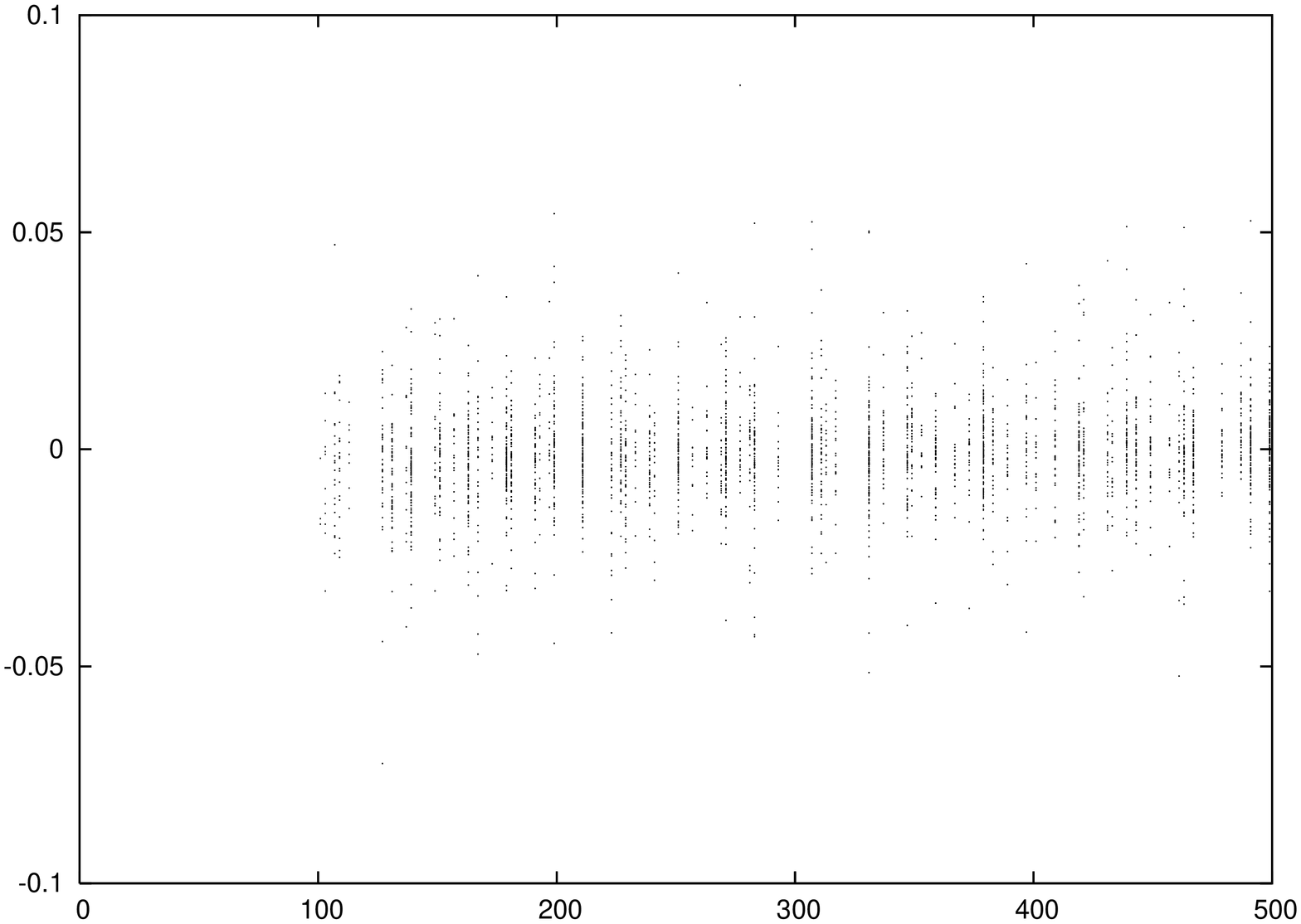,width=1.25in,angle=-90}
            \psfig{figure=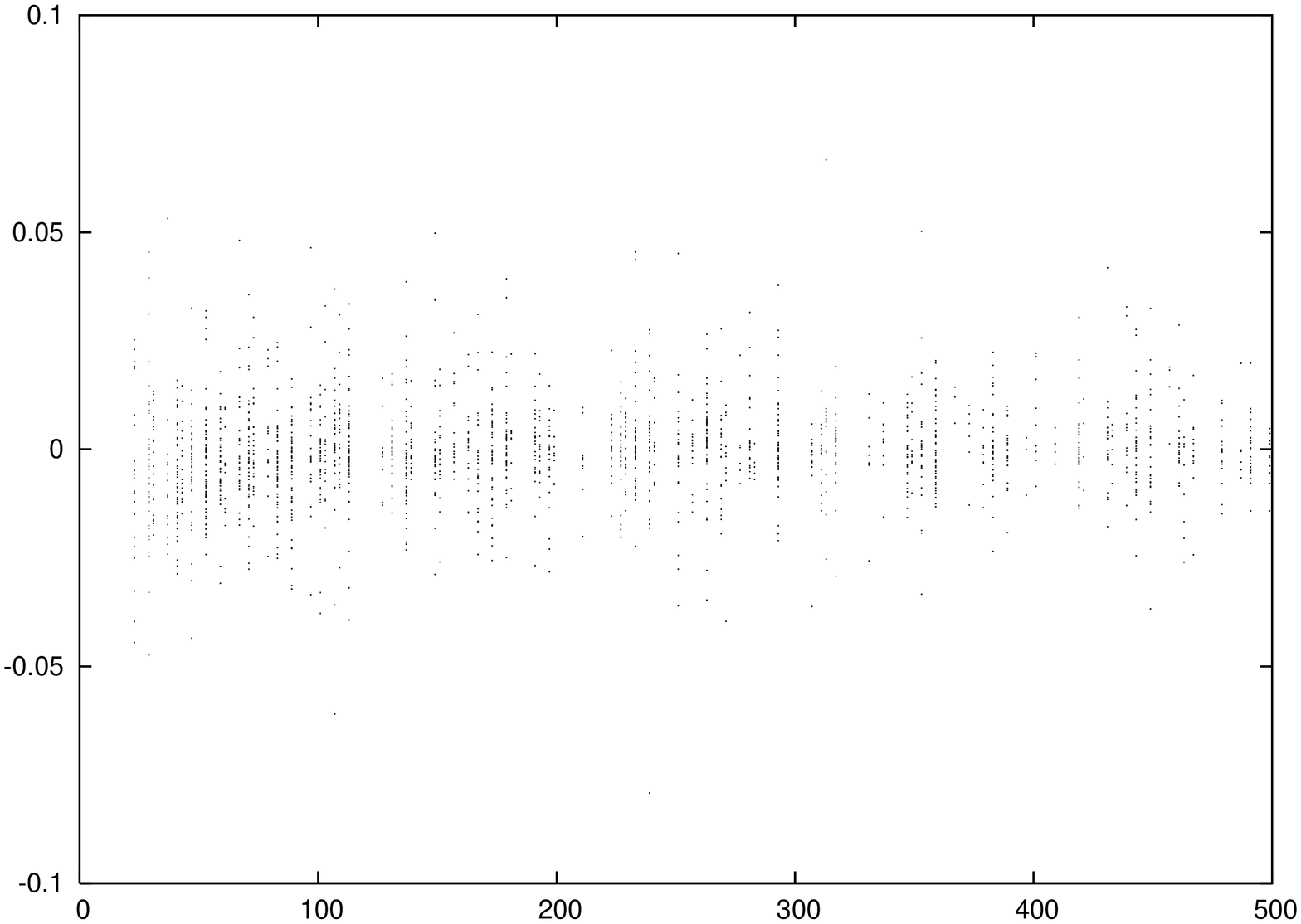,width=1.25in,angle=-90}
            \psfig{figure=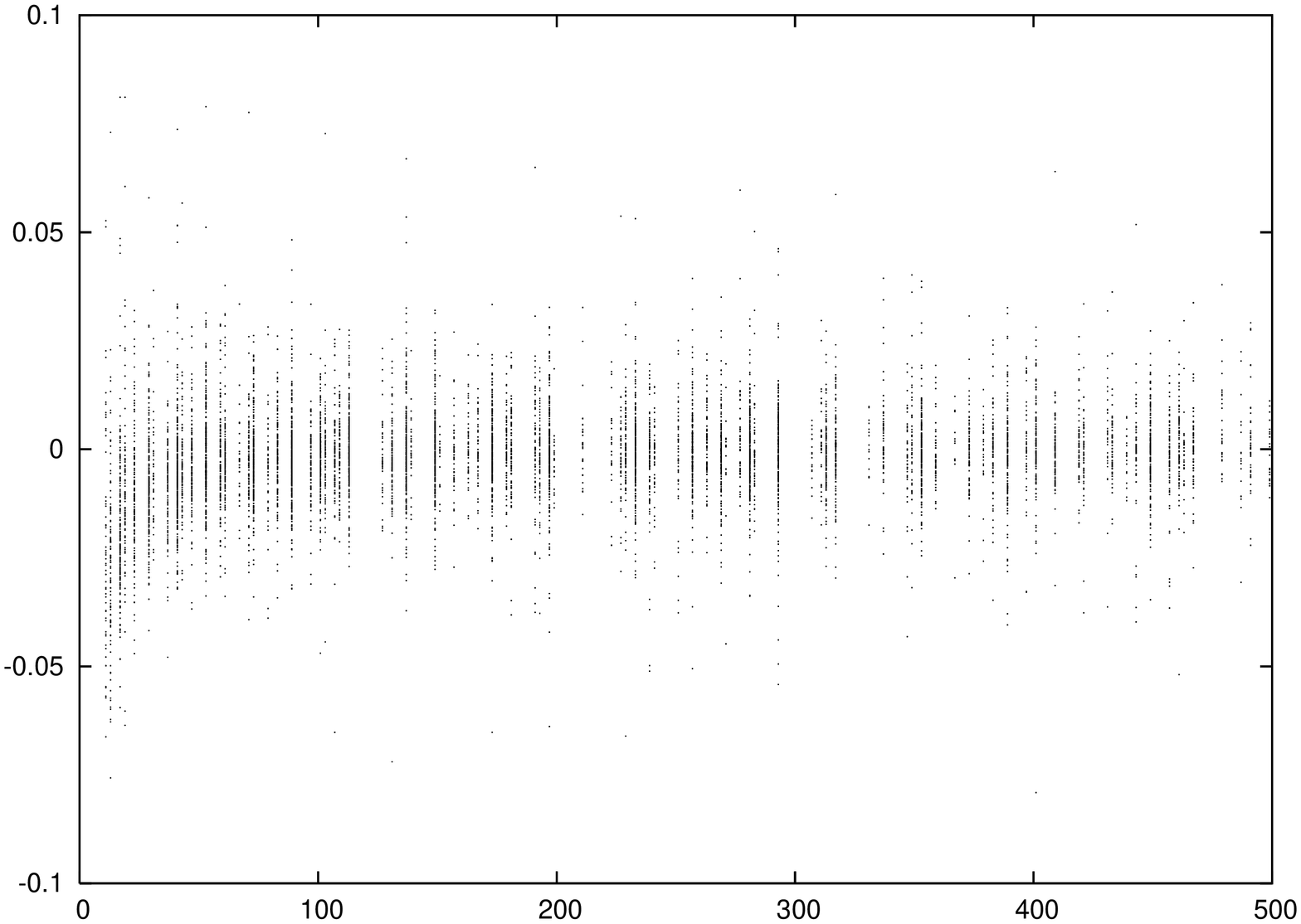,width=1.25in,angle=-90}
    }
    \centerline{
            \psfig{figure=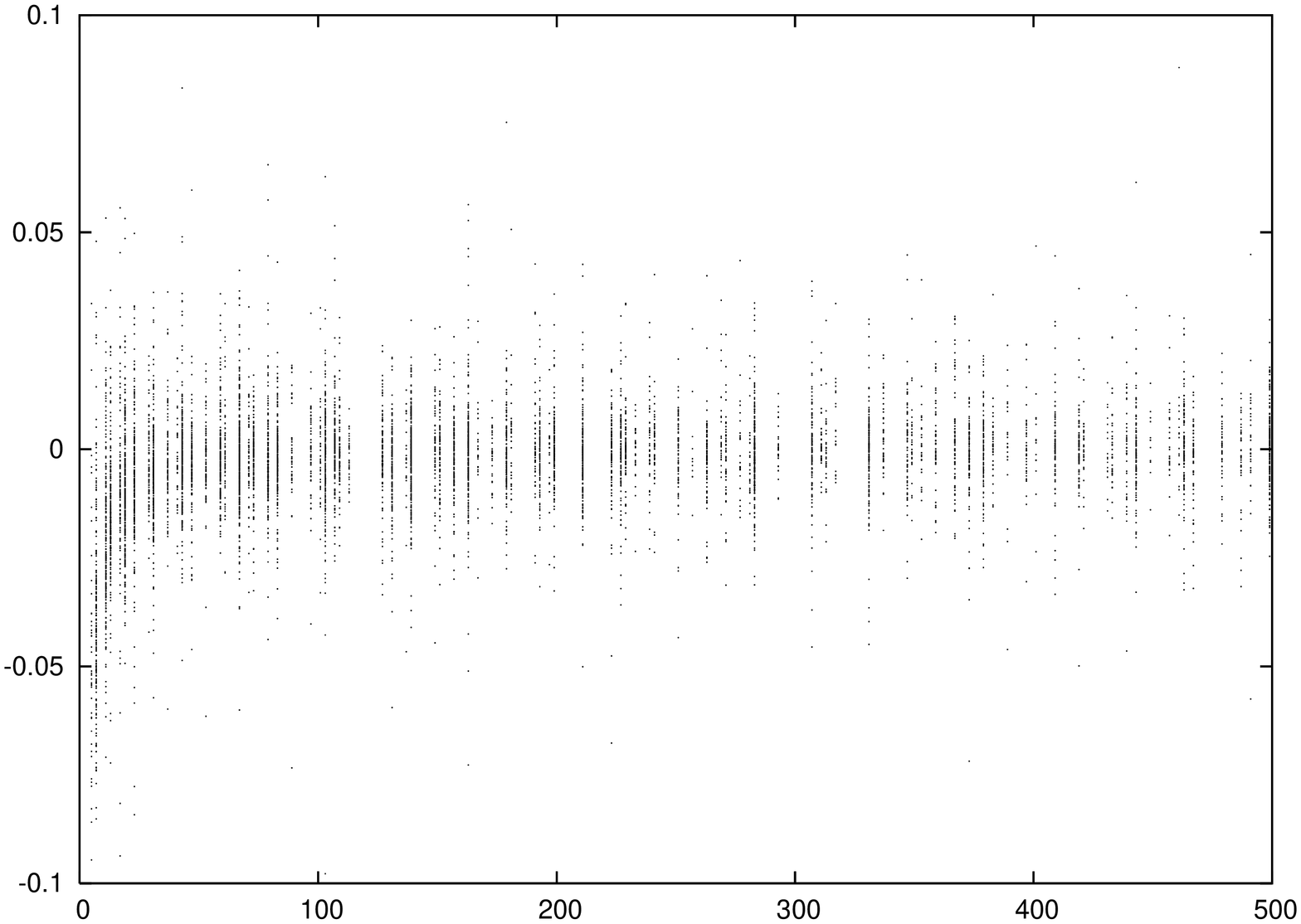,width=1.25in,angle=-90}
            \psfig{figure=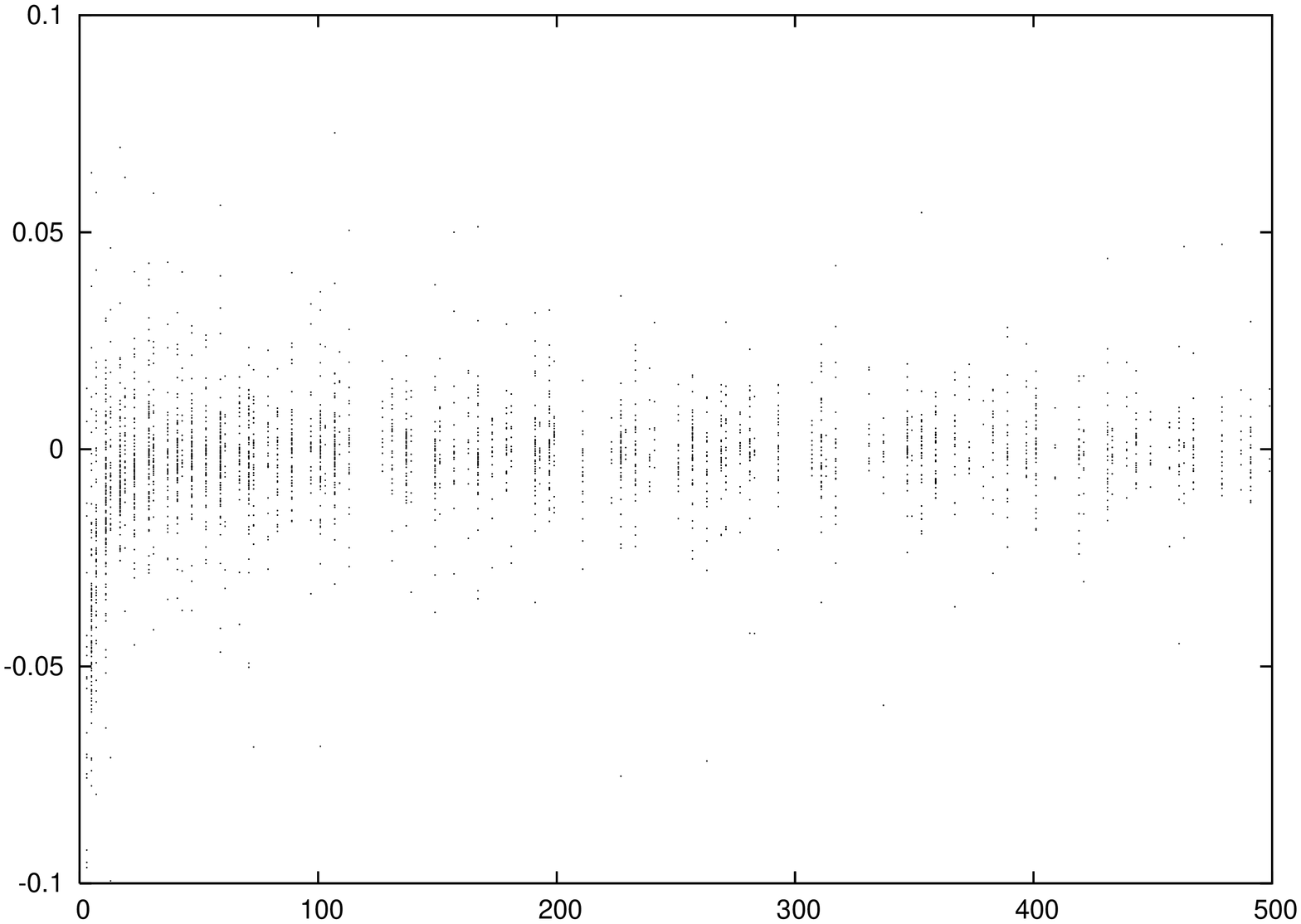,width=1.25in,angle=-90}
            \psfig{figure=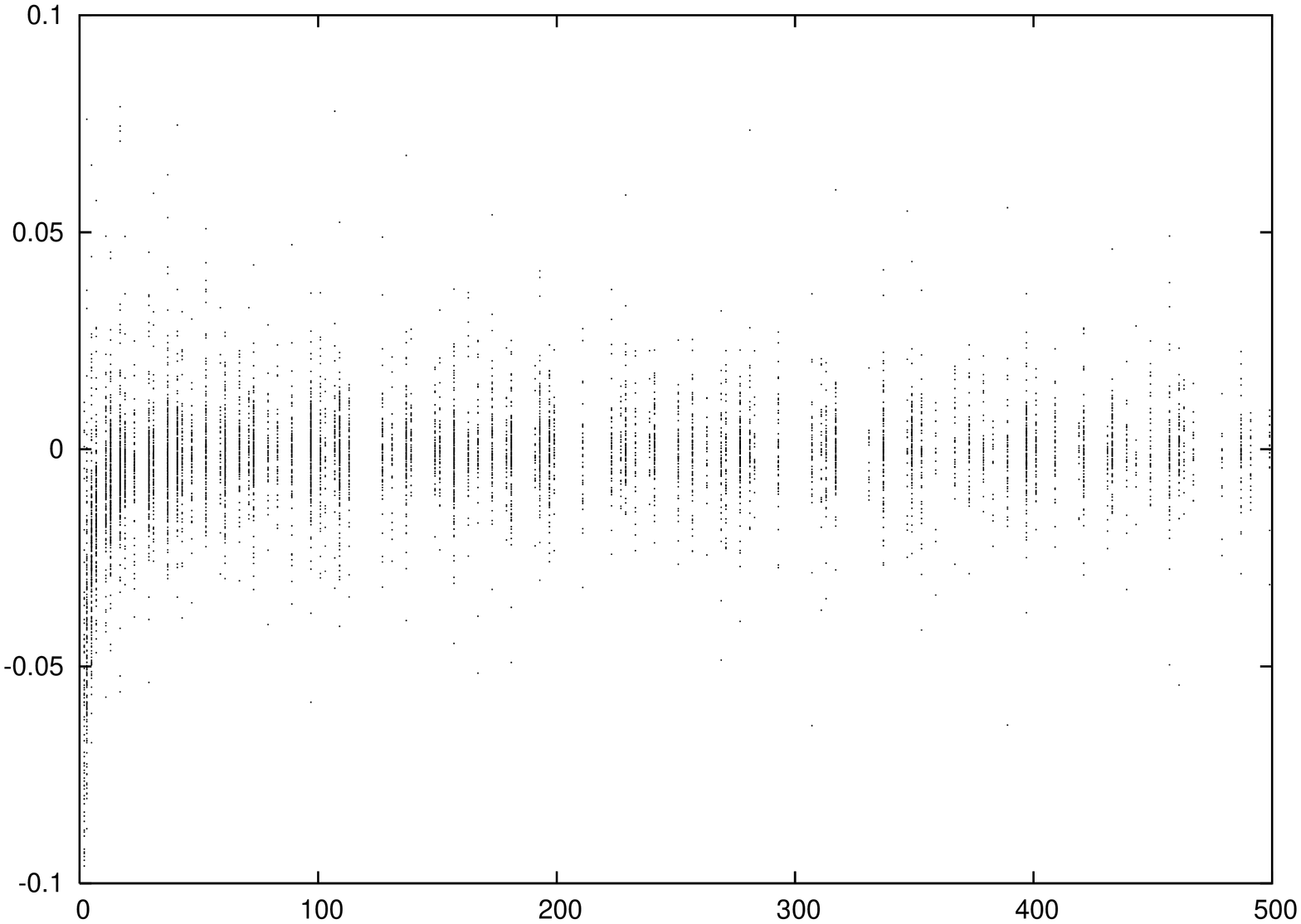,width=1.25in,angle=-90}
    }
    \centerline{
            \psfig{figure=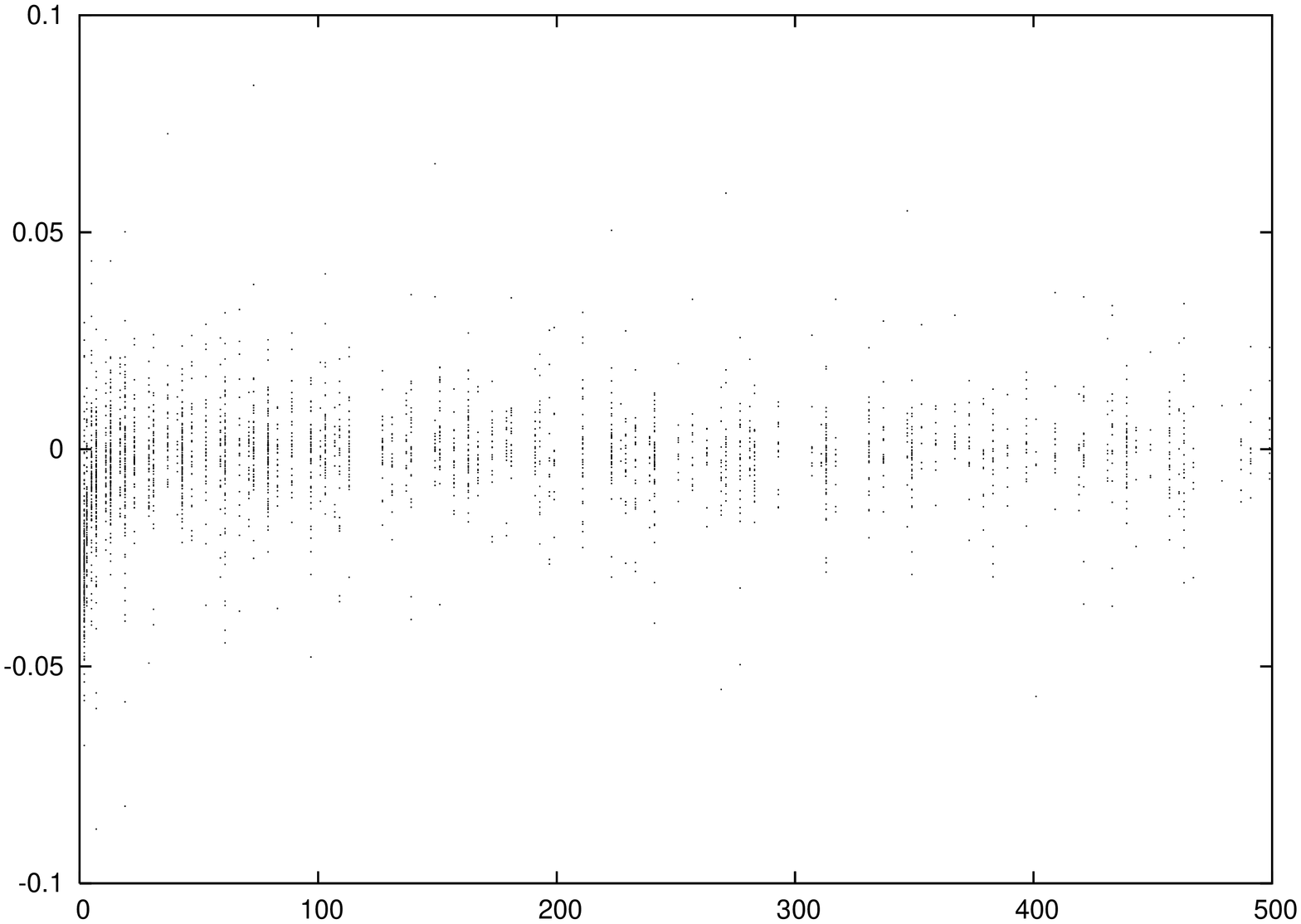,width=1.25in,angle=-90}
            \psfig{figure=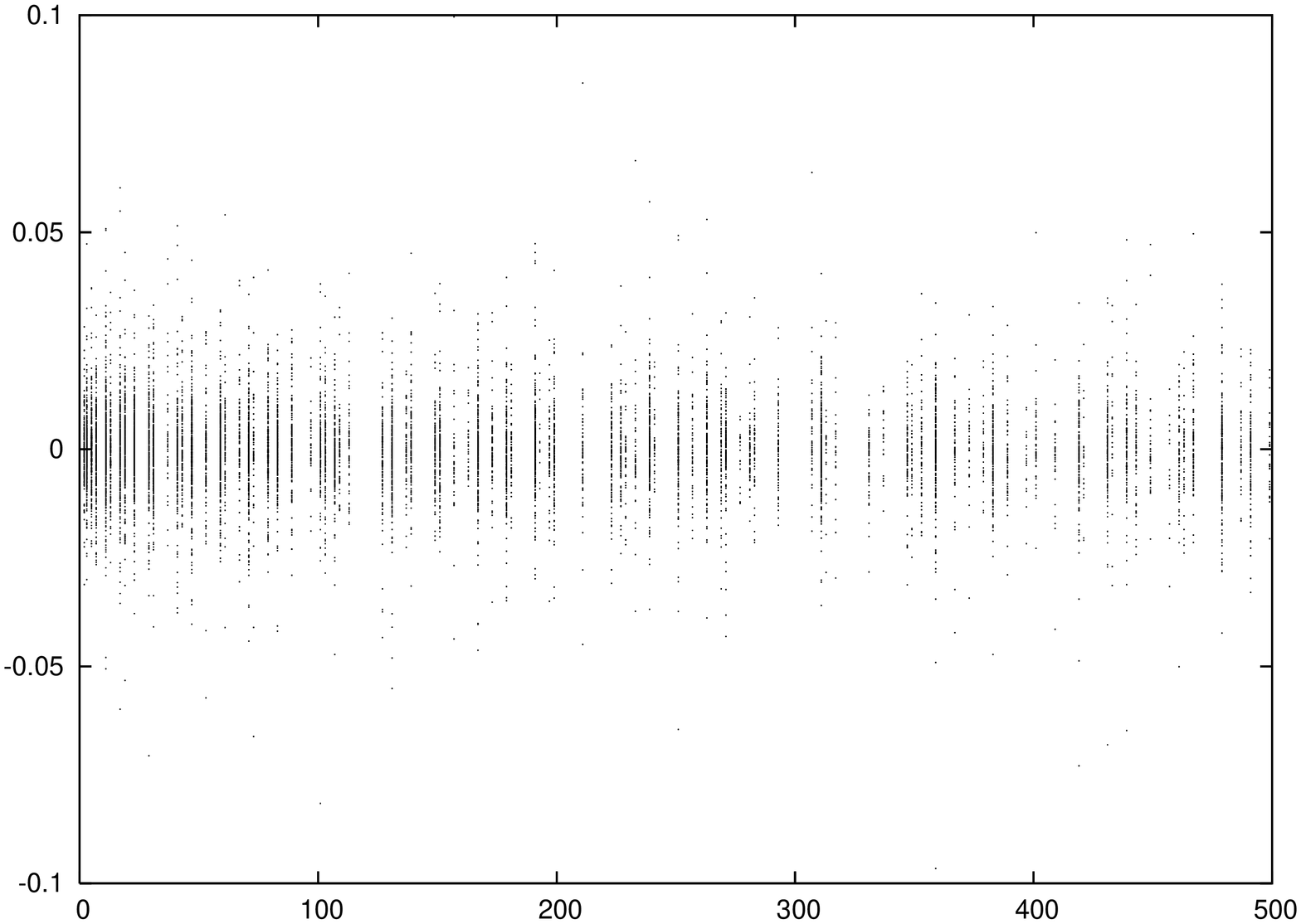,width=1.25in,angle=-90}
            \psfig{figure=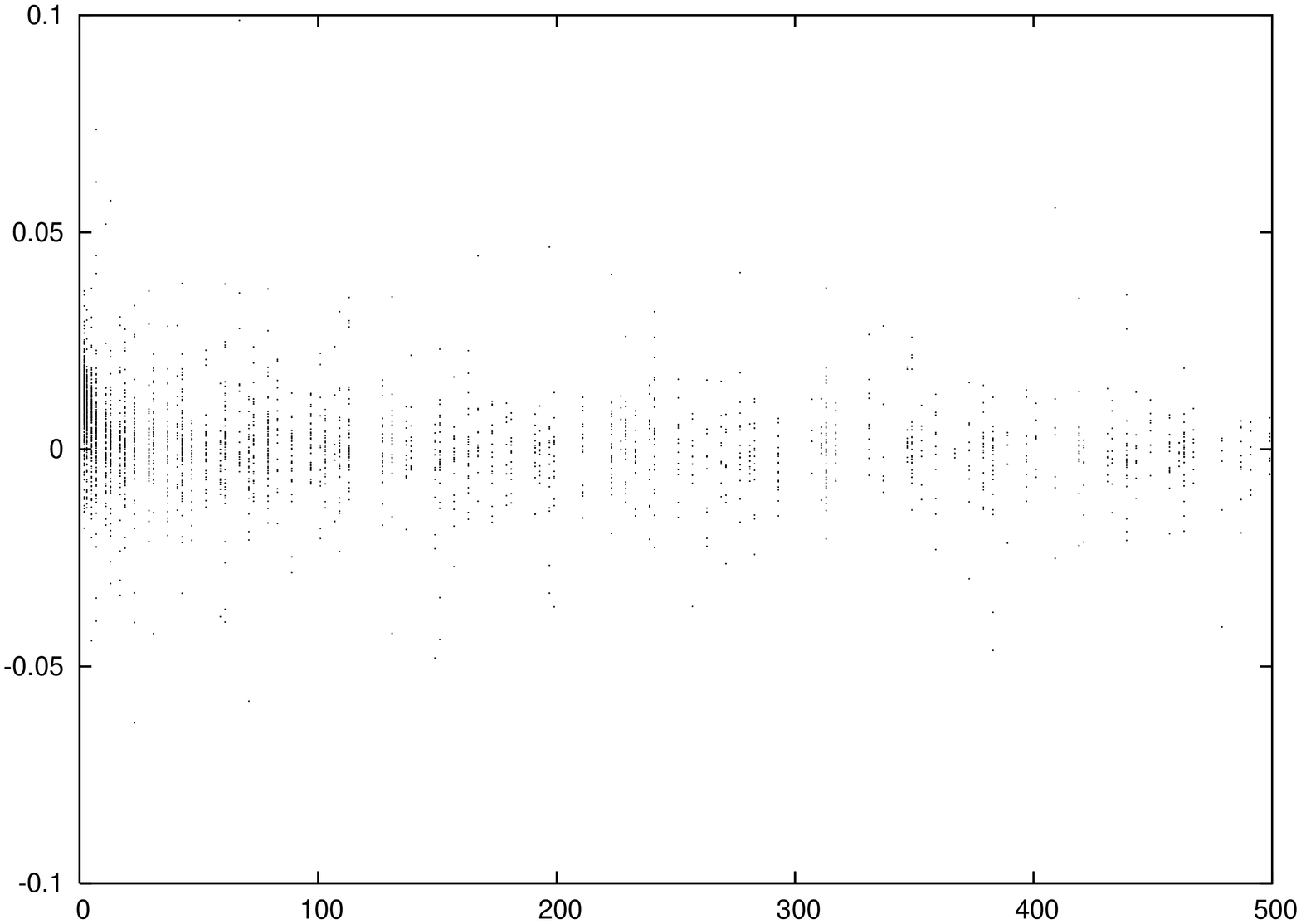,width=1.25in,angle=-90}
    }
    \centerline{
            \psfig{figure=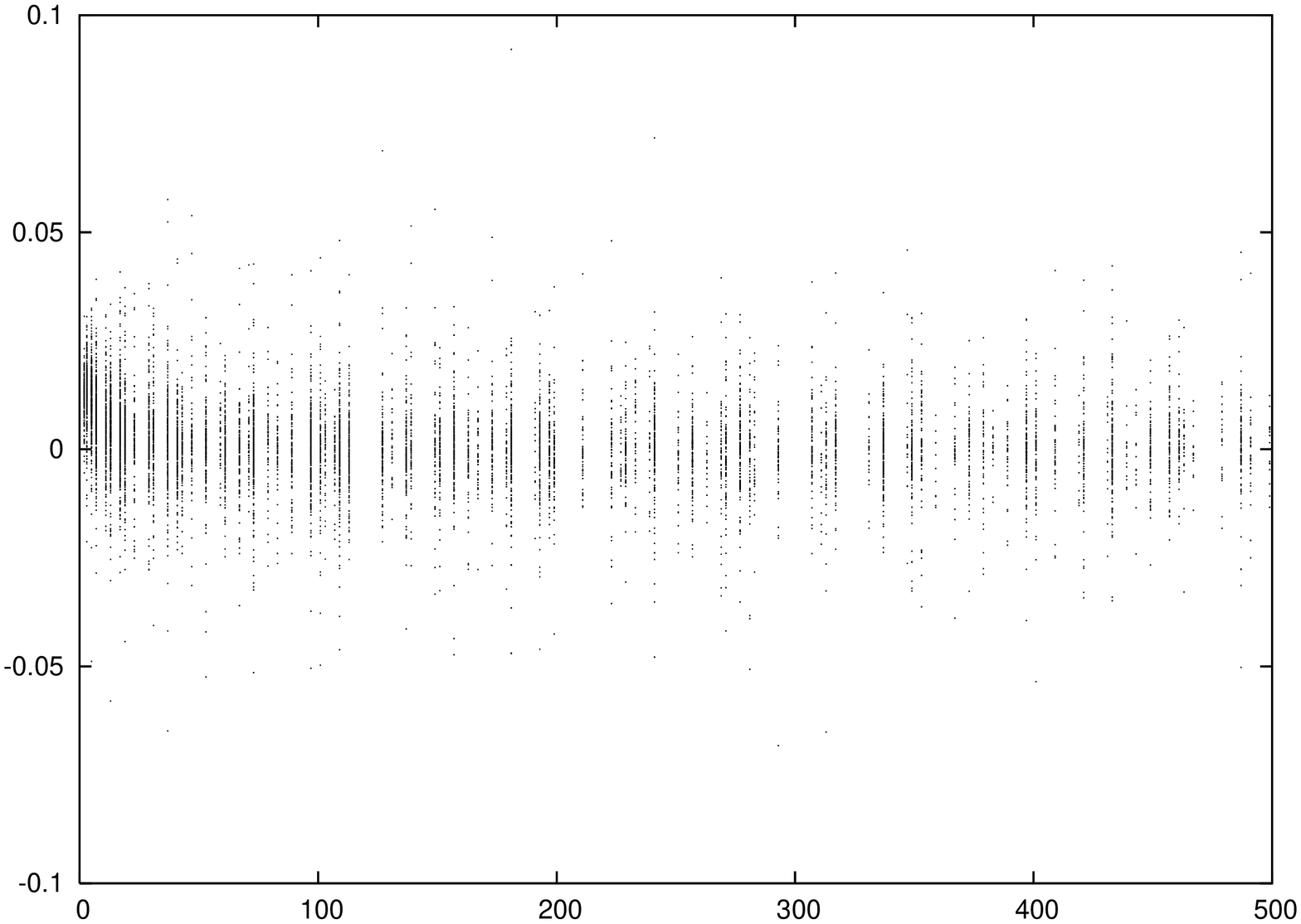,width=1.25in,angle=-90}
            \psfig{figure=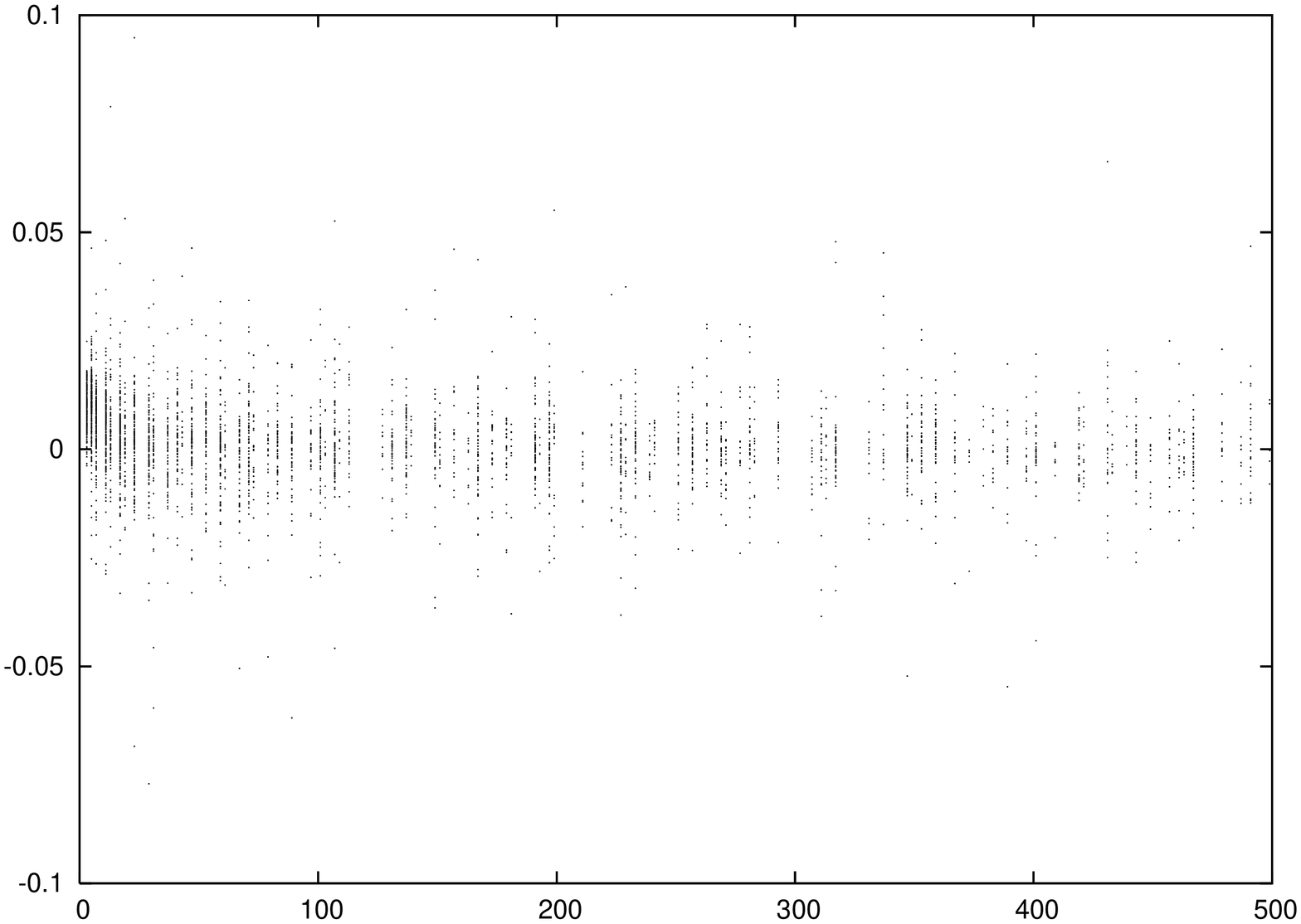,width=1.25in,angle=-90}
            \psfig{figure=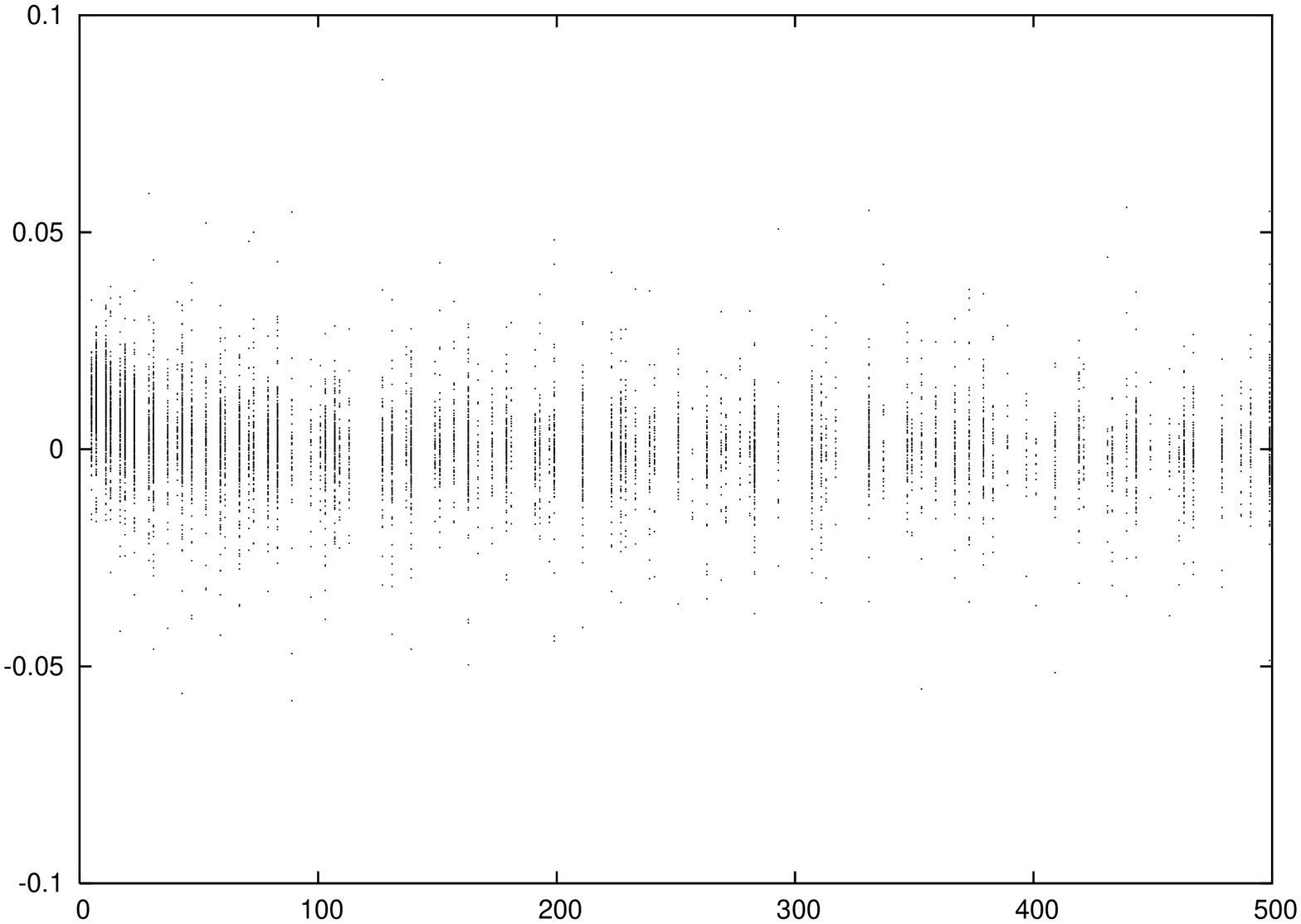,width=1.25in,angle=-90}
    }
    \centerline{
            \psfig{figure=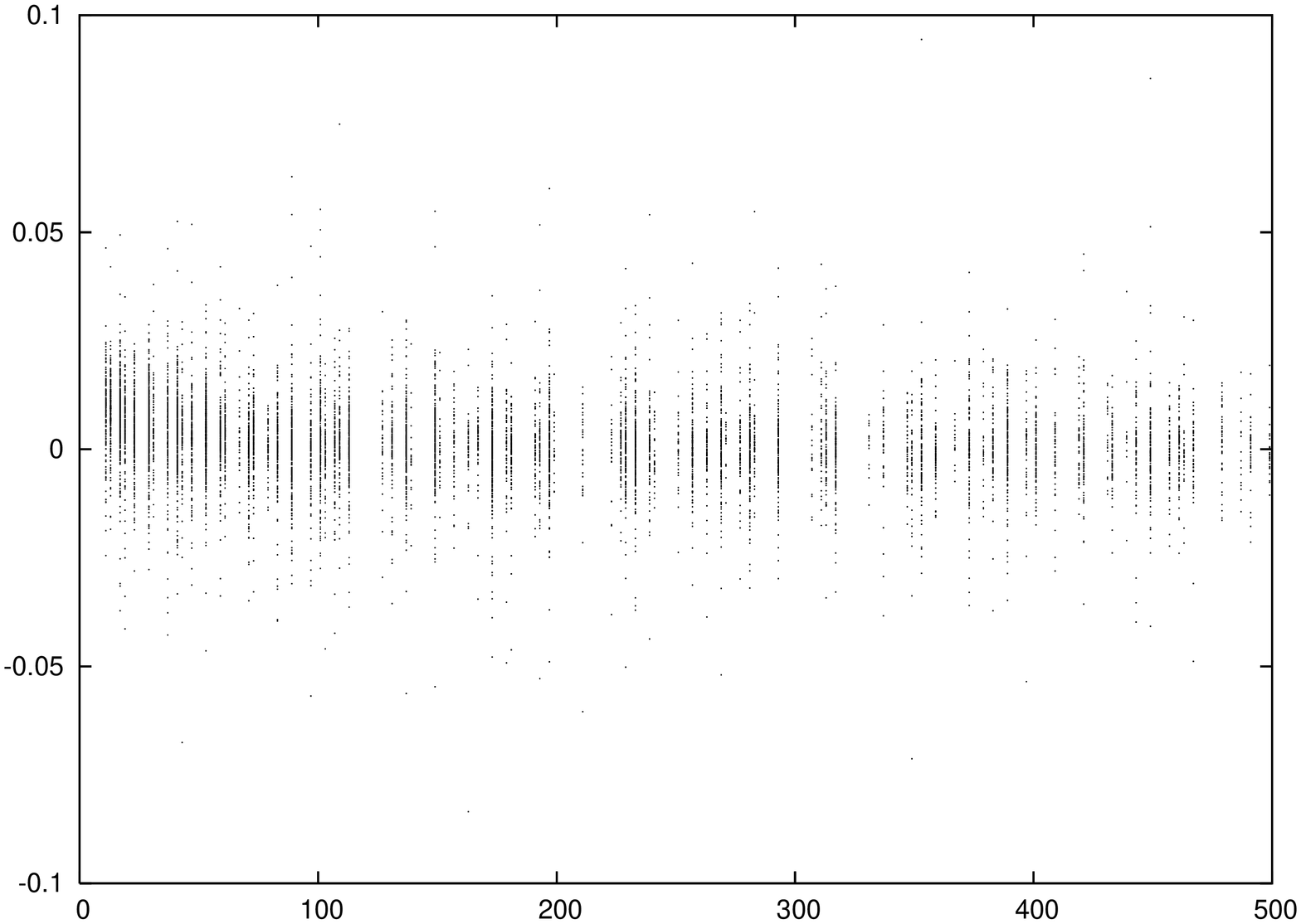,width=1.25in,angle=-90}
            \psfig{figure=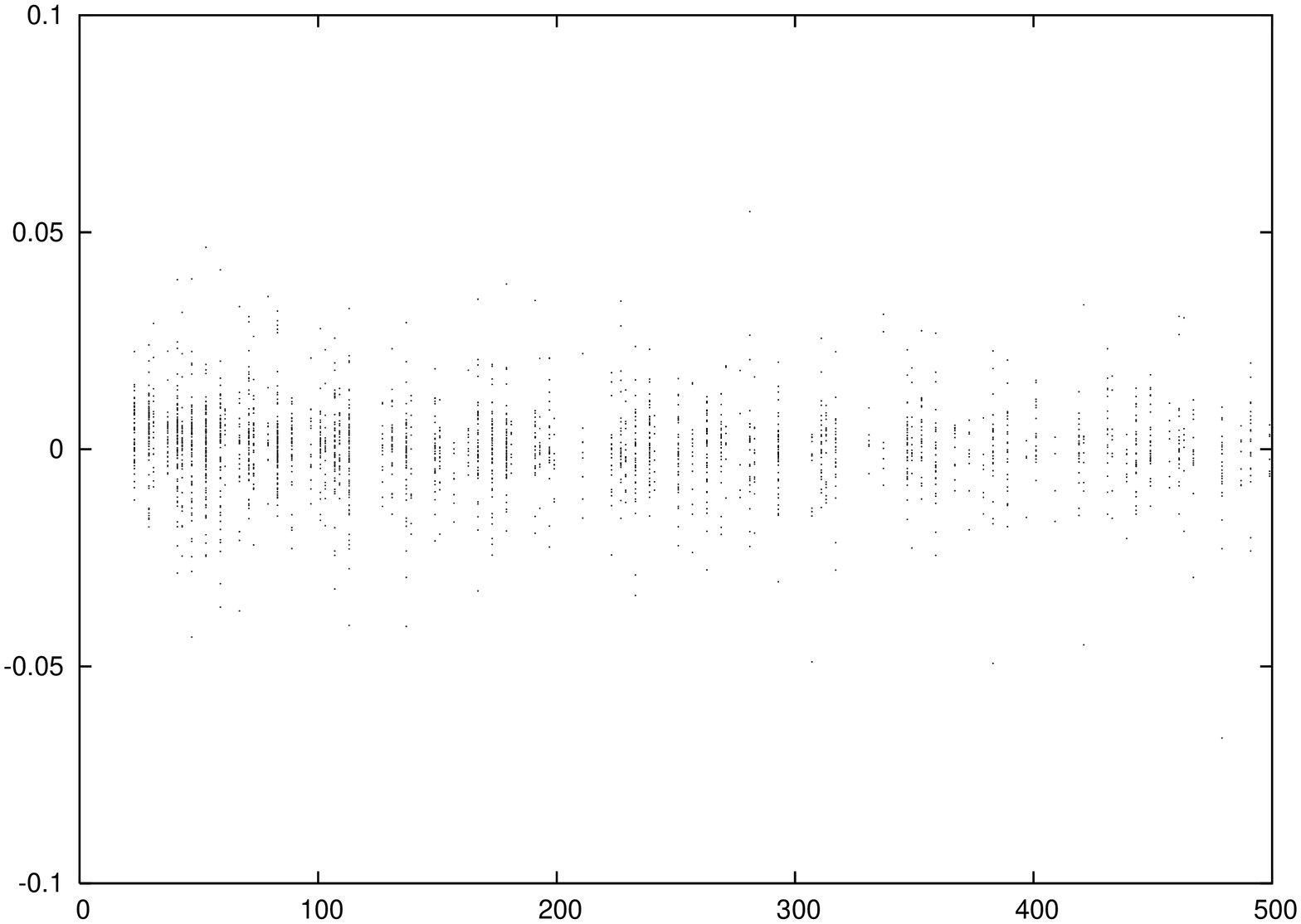,width=1.25in,angle=-90}
            \psfig{figure=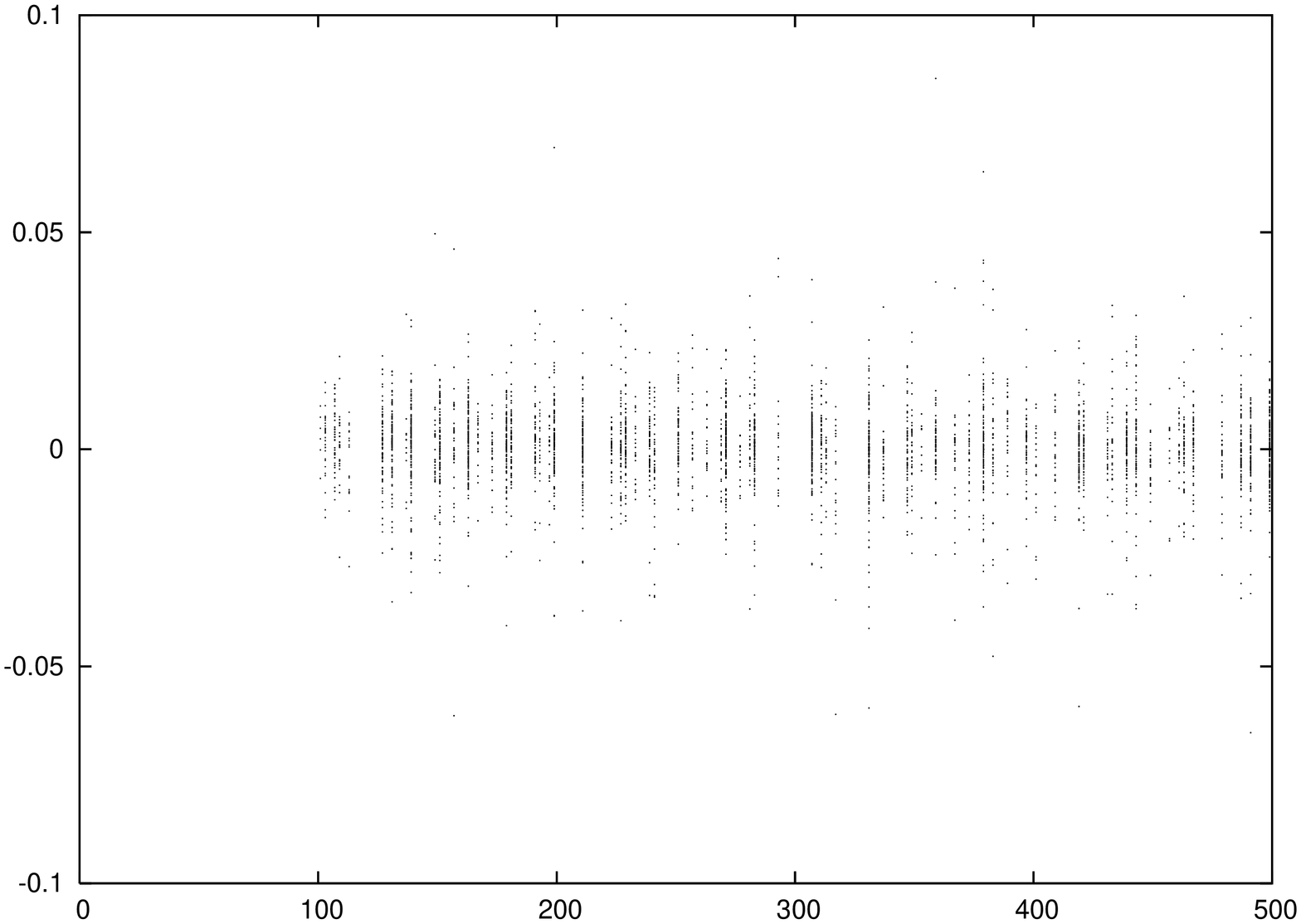,width=1.25in,angle=-90}
    }
    \caption
    {Left to right, top to bottom: $n=-20,-9,-6,-4,-3,-2,-1,0,1,2,3,4,6,9,20$
     Values of (\ref{eq:reality minus conj2}), with $X=10^8$,
     $2\leq q <500$, for the subset of our elliptic curves
     satisfying $a_q=n$.
    }
    \label{fig:seq 2nd}
\end{figure}

Finally, a sequence of plots shows the dependence of the remainder
term in the first and second order approximations on $q$ and $a_q$. Given an integer
$n$, we display, in Figure~\ref{fig:seq 1st}
$q$ v.s. $R^{\pm}_q(10^8)-R_q$ for the subset of our elliptic curves
satisfying $a_q=n$. For each of $n=-20,-9,-6,-4,-3,-2,-1,0,1,2,3,4,6,9,20$ there is one plot.
Figure~\ref{fig:seq 2nd} does the same but for the values given by
equation~(\ref{eq:reality minus conj2}).

We notice several things. Overall, the plots in the Figure~\ref{fig:seq 2nd}
are more symmetric about the horizontal axis reflecting a tighter
fit by including the second order term. For smaller $q$ however, incorporating
the second order term leads to a correction that tends to overshoot. Compare
for example the fourth plot in Figures~\ref{fig:seq 1st} and~\ref{fig:seq 2nd}.
Presumably, the third and further order terms, while of size $O(\log(X)^{-2})$
can have relatively large constants for smaller $q$ requiring one
to take $X$ larger than $10^8$ to see an improvement from the second order term.

This is also reflected in Tables~\ref{tab:11A}--~\ref{tab:307A} which lists
for two elliptic curves and the sets $S^+(10^8)$ and $S^-(10^8)$ the numeric values
of~(\ref{eq:reality minus conj1}) and~(\ref{eq:reality minus conj2}) for
$q \leq 179$.

\begin{table}[h!tb]\begin{small}
\centerline{
\begin{tabular}{|c|c||c|c||c|c|}\hline
$q$ & $a_q$ &
(\ref{eq:reality minus conj1}), $R^-$ case & (\ref{eq:reality minus conj2}), $R^-$ case&
(\ref{eq:reality minus conj1}), $R^+$ case& (\ref{eq:reality minus conj2}),  $R^+$ case\cr
\hline
2 & -2 & -0.0770803072 & -0.1058493733 & -0.0586746787 & -0.0877402111 \cr
3 & -1 & -0.0226715635 & -0.0314020531 & -0.0112745015 & -0.0200944948 \cr
5 & 1 & 0.0039386614 & 0.0110670332 & 0.0036670414 & 0.0108679937 \cr
7 & -2 & -0.0086677613 & -0.0320476479 & 0.0122162834 & -0.0114052128 \cr
13 & 4 & -0.0117312471 & 0.0114581936 & -0.0109800729 & 0.0124435613 \cr
17 & -2 & 0.0068671146 & -0.0078374991 & 0.0156420190 & 0.0007858160 \cr
19 & 0 & 0.0018786796 & 0.0018786796 & 0.0017548761 & 0.0017548761 \cr
23 & -1 & 0.0065085545 & 0.0007253864 & 0.0087254527 & 0.0028829043 \cr
29 & 0 & 0.0015867409 & 0.0015867409 & 0.0024574134 & 0.0024574134 \cr
31 & 7 & -0.0203976628 & 0.0065021478 & -0.0212844047 & 0.0058867043 \cr
37 & 3 & -0.0076213530 & 0.0038881303 & -0.0081586993 & 0.0034679279 \cr
41 & -8 & 0.0293718254 & -0.0104233512 & 0.0370003139 & -0.0032097869 \cr
43 & -6 & 0.0200767559 & -0.0066399665 & 0.0230632720 & -0.0039304770 \cr
47 & 8 & -0.0166158276 & 0.0077120067 & -0.0181946828 & 0.0063789626 \cr
53 & -6 & 0.0175200151 & -0.0048911726 & 0.0194053316 & -0.0032378110 \cr
59 & 5 & -0.0095451504 & 0.0043844494 & -0.0127090647 & 0.0013621363 \cr
61 & 12 & -0.0229944549 & 0.0068341556 & -0.0279181705 & 0.0022108579 \cr
67 & -7 & 0.0114509369 & -0.0104875891 & 0.0227073168 & 0.0005417642 \cr
71 & -3 & 0.0078736247 & -0.0004772247 & 0.0051206275 & -0.0033160932 \cr
73 & 4 & -0.0037492048 & 0.0060879152 & -0.0119406010 & -0.0020032563 \cr
79 & -10 & 0.0300180540 & 0.0013488112 & 0.0296738495 & 0.0007070253 \cr
83 & -6 & 0.0142507227 & -0.0012053860 & 0.0124985709 & -0.0031170117 \cr
89 & 15 & -0.0230738419 & 0.0057929377 & -0.0246777538 & 0.0044799769 \cr
97 & -7 & 0.0105905604 & -0.0054712607 & 0.0154867447 & -0.0007408496 \cr
101 & 2 & -0.0037100582 & 0.0002953972 & -0.0044847165 & -0.0004383257 \cr
103 & -16 & 0.0324024693 & -0.0068711726 & 0.0357260869 & -0.0039571170 \cr
107 & 18 & -0.0228240764 & 0.0073200274 & -0.0245602341 & 0.0058874808 \cr
109 & 10 & -0.0097574184 & 0.0078543625 & -0.0133419792 & 0.0044484844 \cr
113 & 9 & -0.0120886539 & 0.0035056429 & -0.0113667336 & 0.0043859550 \cr
127 & 8 & -0.0093873089 & 0.0034881040 & -0.0081483592 & 0.0048580252 \cr
131 & -18 & 0.0320681832 & -0.0038139100 & 0.0371594888 & 0.0009037228 \cr
137 & -7 & 0.0117897817 & -0.0002445226 & 0.0086451554 & -0.0035131214 \cr
139 & 10 & -0.0148514126 & 0.0000259176 & -0.0112784046 & 0.0037500975 \cr
149 & -10 & 0.0140952751 & -0.0023344544 & 0.0172405748 & 0.0006412396 \cr
151 & 2 & -0.0041170706 & -0.0011557351 & -0.0070016068 & -0.0040099902 \cr
157 & -7 & 0.0108322334 & 0.0000925632 & 0.0097641977 & -0.0010860401 \cr
163 & 4 & -0.0014750980 & 0.0040361356 & -0.0066512858 & -0.0010837710 \cr
167 & -12 & 0.0171132732 & -0.0010302790 & 0.0222297420 & 0.0038987403 \cr
173 & -6 & 0.0054181738 & -0.0030119338 & 0.0036566390 & -0.0048601622 \cr
179 & -15 & 0.0177416502 & -0.0040658274 & 0.0261766468 & 0.0041434818 \cr
\hline
\end{tabular}
} \caption{The values of $R^{\pm}_q(10^8)$, for the elliptic curve
$11_A$ of conductor 11 given by $y^2+y=x^3-x^2-10x-20$, compared
to the conjectured first order approximation~(\ref{eq:reality
minus conj1}) and second order approximation~(\ref{eq:reality
minus conj2}).} \label{tab:11A}\end{small}\end{table}

\begin{table}[h!tb]\begin{small}
\centerline{
\begin{tabular}{|c|c||c|c||c|c|}\hline
$q$ & $a_q$ &
(\ref{eq:reality minus conj1}), $R^-$ case & (\ref{eq:reality minus conj2}), $R^-$ case&
(\ref{eq:reality minus conj1}), $R^+$ case& (\ref{eq:reality minus conj2}),  $R^+$ case\cr
\hline
2 & 0 & 0.0001964177 & 0.0001964177 & 0.0025336244 & 0.0025336244 \cr
3 & 0 & -0.0007380207 & -0.0007380207 & -0.0025236647 & -0.0025236647 \cr
5 & 4 & -0.0128879806 & 0.0109510354 & -0.0166316058 & 0.0072258354 \cr
7 & 0 & -0.0048614428 & -0.0048614428 & -0.0014203548 & -0.0014203548 \cr
11 & 3 & -0.0076239866 & 0.0095824910 & -0.0101221542 & 0.0070977143 \cr
13 & 6 & -0.0212338218 & 0.0089386380 & -0.0276990384 & 0.0024967032 \cr
17 & -1 & 0.0033655021 & -0.0029005302 & 0.0086797465 & 0.0024087738 \cr
19 & -1 & 0.0055745934 & -0.0003223680 & 0.0020465484 & -0.0038550619 \cr
23 & -2 & 0.0074744917 & -0.0036255406 & 0.0079256468 & -0.0031831583 \cr
29 & 0 & 0.0004190042 & 0.0004190042 & -0.0010879108 & -0.0010879108 \cr
31 & 4 & -0.0108662407 & 0.0041956843 & -0.0096223973 & 0.0054512748 \cr
37 & 3 & -0.0067227670 & 0.0037940655 & -0.0162107316 & -0.0056856756 \cr
41 & 5 & -0.0109118090 & 0.0049138186 & -0.0164777387 & -0.0006397725 \cr
43 & -10 & 0.0406071465 & -0.0060036473 & 0.0409949651 & -0.0056532348 \cr
47 & -6 & 0.0284021024 & 0.0057897746 & 0.0209827487 & -0.0016475471 \cr
53 & -10 & 0.0361234610 & -0.0017568821 & 0.0423409405 & 0.0044303004 \cr
59 & 4 & -0.0054935724 & 0.0048495607 & -0.0148985734 & -0.0045473511 \cr
61 & -8 & 0.0227634479 & -0.0025651538 & 0.0253866588 & 0.0000379053 \cr
67 & -8 & 0.0217284008 & -0.0016029354 & 0.0249365465 & 0.0015866634 \cr
71 & -15 & 0.0398795640 & -0.0080932079 & 0.0531538377 & 0.0051425339 \cr
73 & 2 & -0.0003657281 & 0.0042519609 & -0.0019954011 & 0.0026259102 \cr
79 & -13 & 0.0270702549 & -0.0087950276 & 0.0328555729 & -0.0030383756 \cr
83 & 5 & -0.0120289758 & -0.0018576129 & -0.0140337206 & -0.0038544019 \cr
89 & 9 & -0.0117002661 & 0.0050406278 & -0.0159501141 & 0.0008038275 \cr
97 & 7 & -0.0121449601 & 0.0003884458 & -0.0126491435 & -0.0001059465 \cr
101 & 10 & -0.0162655200 & 0.0006799944 & -0.0166873803 & 0.0002713400 \cr
103 & 11 & -0.0154514081 & 0.0027879315 & -0.0155096044 & 0.0027439391 \cr
107 & -15 & 0.0298791131 & -0.0020491232 & 0.0346054275 & 0.0026517125 \cr
109 & -7 & 0.0131301691 & -0.0001660211 & 0.0138662913 & 0.0005595788 \cr
113 & 14 & -0.0219346950 & -0.0006197951 & -0.0199798581 & 0.0013516122 \cr
127 & 17 & -0.0231978866 & 0.0002623636 & -0.0235951007 & -0.0001166344 \cr
131 & -6 & 0.0075864820 & -0.0020988314 & 0.0132703492 & 0.0035773840 \cr
137 & -6 & 0.0049307893 & -0.0044030816 & 0.0085067787 & -0.0008344650 \cr
139 & 14 & -0.0179638452 & 0.0005718014 & -0.0220919830 & -0.0035419036 \cr
149 & 19 & -0.0184587534 & 0.0048182811 & -0.0206858659 & 0.0026092450 \cr
151 & -14 & 0.0157624561 & -0.0057016072 & 0.0217272592 & 0.0002461455 \cr
157 & -14 & 0.0258394912 & 0.0051129620 & 0.0236949594 & 0.0029519710 \cr
163 & -8 & 0.0115026664 & 0.0005637031 & 0.0044174198 & -0.0065301910 \cr
167 & 21 & -0.0224356707 & 0.0011192167 & -0.0284090909 & -0.0048359139 \cr
173 & -6 & 0.0056158047 & -0.0020804090 & 0.0044893748 & -0.0032129134 \cr
179 & 0 & 0.0018844544 & 0.0018844544 & -0.0007004350 & -0.0007004350 \cr
\hline
\end{tabular}
} \caption{The values of $R^{\pm}_q(10^8)$, for the elliptic curve
$307_A$ of conductor 307 given by $y^2 +y= x^3-x-9$, compared to
the conjectured first order approximation~(\ref{eq:reality minus
conj1}) and second order approximation~(\ref{eq:reality minus
conj2}).} \label{tab:307A}\end{small}\end{table}

\vspace{.1in}

{\noindent \bf Acknowledgements}

We wish to thank the Newton Institute in Cambridge where some of this research was carried out.
We also thank Gonzalo Tornaria and Fernando Rodriguez-Villegas who supplied us with
a table of weight three halves forms that were used to compute the $L$-values studied
in this paper.

\begin{flushleft}
\mbox{}\\
\mbox{}\\
J.B. Conrey\\American Institute of Mathematics\\ 360 Portage
Avenue\\Palo Alto, CA 94306\\USA\\
\mbox{}\\
School of Mathematics \\ University of Bristol\\ Bristol BS8 1TW\\
UK\\
\mbox{}\\
\mbox{}\\
A. Pocharel\\
Department of Mathematics \\ Princeton University\\ Princeton, NJ
08544 \\ USA \\
\mbox{}\\
\mbox{}\\
M.O. Rubinstein\\
Pure Mathematics \\ University of Waterloo\\ 200 University
Ave W\\Waterloo, ON, Canada\\N2L 3G1\\
\mbox{}\\
\mbox{}\\
M. Watkins\\
School of Mathematics \\ University of Bristol\\ Bristol BS8 1TW\\
UK
\end{flushleft}

\end{document}